\documentclass[twocolumn,10pt]{article}

\usepackage{newtxtext,newtxmath}   
\usepackage{microtype}
\usepackage{booktabs}
\usepackage{longtable}
\usepackage{multicol}

\usepackage[a4paper,margin=1in]{geometry}
\usepackage{setspace}
\usepackage{fancyhdr}
\usepackage{titlesec}
  \titleformat{\section}
    {\normalfont\bfseries\fontsize{12}{15}}{\thesection}{1em}{}
  \titlespacing*{\section}{0pt}{10pt}{5pt}
  \titleformat{\subsection}
    {\normalfont\bfseries\fontsize{11}{13}}{\thesubsection}{1em}{}
  \titlespacing*{\subsection}{0pt}{8pt}{4pt}

\usepackage{graphicx}
  \graphicspath{{./images/}{./ANT_Problems/}{./}}
  \DeclareGraphicsExtensions{.png,.pdf,.jpg}

\usepackage{amsmath}

\usepackage{hyperref}
\hypersetup{colorlinks=true,linkcolor=blue,citecolor=blue,urlcolor=blue}

\usepackage{xcolor}

\usepackage{booktabs,tabularx}

\usepackage{listings}

\usepackage{algorithm,algpseudocode}
\usepackage{mdframed}
\newmdenv[
  backgroundcolor=gray!10,
  linecolor=gray!60,
  linewidth=1pt,
  roundcorner=4pt,
  innertopmargin=5pt
]{block}
\usepackage{enumitem}
\usepackage{cite}
\usepackage{caption}
\title{Computational Experiments in Number Theory}
\author{Ali Saraeb}
\date{April 2025}

\begin{document}

\maketitle
\fancypagestyle{firstpage}{%
  \fancyhf{}                    
  \renewcommand{\headrulewidth}{0pt}
  \renewcommand{\footrulewidth}{0pt}
  \fancyfoot[C]{%
    \footnotesize
    Department of Mathematics, Ohio State University\\
    \texttt{saraeb.1@osu.edu , ali.y.saraeb@gmail.com}%
  }%
}
\thispagestyle{firstpage}

\begin{abstract}
This paper presents two concrete applications of Artificial Intelligence to algorithmic and analytic number theory.

Recent benchmarks of large language models have mainly focused on general mathematics problems and automated theorem-proving. In the first part of this paper, we take a different approach, and we focus on a more specialized domain: we evaluate the performance of the state-of-the-art open source large language model \textit{Qwen2.5‑Math‑7B‑Instruct} on algorithmic and computational tasks in algorithmic number theory. On a benchmark of thirty \textit{algorithmic problems} and thirty \textit{computational questions} taken from classical number-theoretic textbooks and Math StackExchange, the model achieves at least \(0.95\) accuracy (relative to the  "true" answer) on every problem/question when given an optimal \textit{non-spoiling} hint. Moreover, for a fixed \textit{hinting strategy}, the mean accuracy peaks at \(0.88\) when averaged over algorithmic problems and at \(0.89\) when averaged over computational questions, indicating the sensitivity of performance to the choice of hinting strategy. Finally, we refer to our manually constructed dataset as the \emph{Hinted Algorithmic Number Theory} (HANT) dataset and make both the dataset and accompanying code publicly available at \url{ https://doi.org/10.5281/zenodo.15293187} \cite{Saraeb2025Zenodo}.

The second part of the paper empirically verifies a folklore conjecture in analytic number theory stating that “the \textit{modulus} \(q\) of a \textit{Dirichlet character} \(\chi\) is uniquely determined by the initial non-trivial zeros \(\{\rho_1,\dots,\rho_k\}\) (for some $k \in \mathbb{N}$) of the corresponding \textit{Dirichlet $L$-function} \(L(s,\chi)\)”. We form the dataset of zeros from the LMFDB database and turn the conjecture's statement into a multiclass classification problem where the feature vectors are engineered from the initial zeros, and the labels are the corresponding moduli. We train an \textit{LightGBM} multiclass classifier to predict the modulus $q$ for $214$ randomly chosen Dirichlet $L$-functions from a vector of statistical features of their initial zeros (moments, finite-difference statistics, FFT magnitudes, etc. ...). The model empirically verifies the conjecture for small $q$ in the sense that it achieves high test accuracy of at least $93.9 \%$, provided that sufficient statistical properties of the zeros are incorporated in the training process. Based on the empirical results, we propose two new conjectures: 
\begin{enumerate}[label=(\roman*),nosep,leftmargin=*,align=left]
  \item There exist hidden statistical patterns in the nontrivial zeros of Dirichlet \(L\)-functions.
  \item There is an underlying statistical connection among the zeros of Dirichlet \(L\)-functions corresponding to characters of the same modulus.
\end{enumerate}
For the second part of the paper, the code and dataset can be found at \url{https://doi.org/10.5281/zenodo.15460772} \cite{Saraeb2025ZenodoAnalytic}
\end{abstract}

\section{Introduction}
\noindent
Large language models (LLMs) have demonstrated strong capabilities in code generation, open-domain question answering, and quantitative reasoning tasks \cite{brown2020language,lewkowycz2022solving}. More recently, researchers have begun to evaluate LLMs such as \textit{GPT-4} \cite{achiam2023gpt}, \textit{PaLM 2}, and \textit{Qwen2.5-Math-7B-Instruct} \cite{yang2024qwen2} on benchmarks ranging from competition-level problem solving to formal theorem proving \cite{hendrycks2021measuring,satpute2024can,agrawal2024give}. However, a few significant gaps remain less explored. First, although graduate or undergraduate students in mathematics routinely rely on hints to prove a theorem or solve a computational problem, the study of “hint-based” prompting seems to be relatively underexplored in the literature. Second, the intersection of AI and algorithmic number theory, which encompasses standard algorithms in mathematics such as integer and ideal factorization, root-extraction in groups, discrete logarithms, primality testing, etc. $\dots$ \cite{cohen2013course}, does not seem to be thoroughly explored. In the first part of this paper, thus, we address these gaps by evaluating the LLM \textit{Qwen2.5‑Math‑7B‑Instruct} on a benchmark of thirty algorithmic problems and thirty computational questions in algorithmic number theory. Specifically, we consider the two index sets
\begin{align*}
& \{(\text{algorithmic problem}_i,\;\text{hint}_j)\}_{1\le i\le 30,\;1\le j\le 9}
\quad\text{and} \\
&\{(\text{computational question}_i,\;\text{hint}_j)\}_{1\le i\le 30,\;1\le j\le 9},
\end{align*}
for a total of \(30\times 9 + 30\times 9  = 540\) prompts. We use the term “algorithmic problems” for tasks in algorithmic number theory whose solution must be an explicit algorithm, and use “computational questions” for those whose answer should be a computational output. To summarize the performance of the model, let \(a_{i,j}\) denote the accuracy of the model on item \(i\) (algorithmic problem or computational question) when given hint \(j\), for \(i=1,\dots,60\) and \(j=1,\dots,9\).  Then the model achieves the following
\[
\max_{1\le j\le9} \frac{1}{60}\sum_{i=1}^{60} a_{i,j} \;\approx\; 0.885,
\]
and more importantly
\[
\forall\,i\in\{1,\dots,60\},\;\exists\,j\in\{1,\dots,9\}\quad\text{s.t}\quad a_{i,j}\ge0.95.
\]
\textbf{Remark 1.1} A comprehensive description of the nine hinting strategies, the complete list of algorithmic problems and computational questions, and the technical details of the model \textit{Qwen2.5-Math-7B-Instruct} used in Part I of the paper are provided in Appendix A.

In the second part of this paper, we explore the use of classical machine‐learning methods— random‐forest classifiers \cite{breiman2001random} and the gradient‐boosted‐tree ensemble LGBMClassifier \cite{ke2017lightgbm}—in analytic number theory \cite{apostol2013introduction,davenport2013multiplicative}. Building on recent work that applies similar techniques to number theoretic problems \cite{shanker2012neural,he2022machine}, we empirically verify a widely believed conjecture in analytic number theory stating that “the modulus \(q\) of a Dirichlet character \(\chi\) is uniquely determined by the initial non-trivial zeros \(\{\rho_1,\dots,\rho_k\}\) (for some $k \in \mathbb{N}$) of the corresponding Dirichlet $L$-function \(L(s,\chi)\)” for small modulus \(q\) \cite{apostol2013introduction, davenport2013multiplicative}. To test this, we treat the conjecture as a multiclass classification problem with the imaginary part of the zeros as input features. We begin by experimenting on a sample of 21 Dirichlet \(L\)-functions as our initial test set. A Random Forest trained on raw zero sequences managed only about \(61.9\%\) accuracy, but after extending each zero vector with engineered features—statistical moments, finite‐difference metrics, FFT magnitudes, and mean prime gap—the same Random Forest achieved a perfect accuracy (\(100\%\)) on that small sample. Motivated by this observation, we then trained on a larger dataset and tested on 214 randomly chosen \(L\)-functions and switched to an \texttt{LGBMClassifier} with similar feature pipeline. On the larger dataset, the classifier reached a high test accuracy of $0.939$ and a validation accuracy of $0.967$, where the small decrease in accuracy is expected due to the noise introduced by the truncation of the zeros’ decimal parts. Beyond providing strong empirical evidence for the conjecture, this methodology introduces a new framework for pattern discovery in sequences of numbers: the jump in accuracy after incorporating engineered statistical properties indicates the existence of underlying statistical patterns in the zeros of Dirichlet $L$-functions, and the fact that the testing process is carried out on Dirichlet $L$-functions different than those used for training purposes indicates the existence of an underlying pattern that different Dirichlet $L$-functions of same modulus $q$ possess. 

\textbf{Remark 1.2} All requisite mathematical preliminaries and formal definitions for Part II of this paper are provided in Appendix B.

\section{Related Work}
Existing evaluations of LLMs on mathematical reasoning employ a variety of benchmark datasets spanning different levels of difficulty. For instance, the MATH dataset \cite{hendrycks2021measuring} contains competition‐style problems from middle‐ and high‐school levels, MiniF2F contains Olympiad‐level questions, and MathEval collects high‐school competition and curriculum‐level tasks. Leading models, including \textit{GPT‐3} \cite{brown2020language}, \textit{GPT‐4}, \textit{Codex}, \textit{LLaMA 3}, \textit{Minerva} \cite{lewkowycz2022solving}, and \textit{Qwen2.5‐Math‐7B‐Instruct} \cite{yang2024qwen2}, achieve state‐of‐the‐art performance on the aforementioned datasets. Two recent studies examine the effect of hint‐based prompting: Agrawal et al.~\cite{agrawal2024give} show that providing a single hint per problem on MATH dataset yields approximately $56\%$ accuracy for Instruct tuned models (with scores in $[30 \%, 70 \%] $ across prompts), while Fu et al.~\cite{fu2024hint} introduce the HSPMATH dataset and report that the hint‐fine‐tuned LLaMA‐7B model attains $64.3 \%$ accuracy, outperforming both GPT‐3.5 and WizardMath‐13B.

Moreover, recent studies at the intersection of machine learning and number theory have employed classical machine learning algorithms—including neural networks and Bayesian classifiers to predict the nontrivial zeros of the Riemann zeta function on the critical line \cite{shanker2012neural}, and empirically verify the Sato–Tate conjecture \cite{he2022machine}.

\section{Datasets}
\noindent
In Part I, we form the \emph{Hinted Algorithmic Number Theory} (HANT) dataset from thirty number theoretic classical algorithmic tasks and their thirty computational counterparts, taken from classical textbooks on the subject \cite{cohen2013course,shoup2009computational} and Math StackExchange (see the list below). The dataset comprises sixty text documents: the first thirty—\emph{algorithmic problems}—each ask for designing an algorithm to solve a classical number-theoretic task; the remaining thirty—\emph{computational questions}—mirror those tasks by requesting an explicit solution for the corresponding problem under the same topic. Moreover, each of the sixty text documents is divided into eleven sections:  
\begin{enumerate}[label=(\arabic*)]
  \item the problem or question statement,  
  \item–(10) nine pedagogical hints (one per hinting strategy as per Appendix $A$),  
  \item the true algorithm or solution.  
\end{enumerate}
The datasets are used to form $60\times9=540$ prompts by combining, for each problem/question, section (1) with each of sections (2)–(10), and we reserve section (11) exclusively as the gold‐standard answer for evaluation purposes. The sixty tasks are named, indexed, and grouped into nine thematic categories as in Table~1 in Appendix~A:  
\begin{itemize}
  \item \emph{GCD \& Linear Congruence}
  \item \emph{Diophantine Equations \& Continued Fractions}
  \item \emph{Integer Representations}
  \item \emph{Modular Arithmetic \& Root Extraction in Groups}
  \item \emph{Exponentiation \& Discrete Logarithm}
  \item \emph{Primality Testing}
  \item \emph{Sieves \& Totient Functions}
  \item \emph{Factorization of Integers \& Ideals}
  \item \emph{Cryptographic Keys \& Elliptic‐Curve Operations}
\end{itemize}

In Part II, we construct two labeled datasets for the classification problem using the first 25 nontrivial truncated zeros of Dirichlet \(L\)-functions as inputs, as retrieved from the LMFDB database \cite{lmfdb}.  Each datum is of the form
\[
q_\chi\;:\;\bigl(\Im(\rho_1),\Im(\rho_2),\dots,\Im(\rho_{25})\bigr),
\]
where \(q_\chi\) is the modulus of the non-principal Dirichlet character \(\chi\) modulo \(q\).  We restrict to moduli \(q\le200\). For each modulus (label)~$q$, we train on fewer than $\phi(q) - 3$ Dirichlet characters, where $\phi$ denotes Euler's totient function. We evaluate performance using a single test feature vector for each $q \leq 100$, and two test feature vectors for each $100 < q \leq 200$. 

\section{Method}
\subsection{Part I: Prompting the LLM with Algorithmic Number-Theory Tasks}

\begin{table}[!ht]
\centering
\small
\begin{tabular}{@{}ll@{}}
\toprule
\textbf{Attribute}           & \textbf{Specification}                                    \\
\midrule
Model name                   & Qwen2.5-Math-7B-Instruct                                  \\
Parameter count              & $7\times10^{9}$                                           \\
Architecture                 & Decoder-only Transformer                                  \\
Maximum context length       & 8\,192 tokens                                             \\
Release date                 & February 2025                                             \\
\bottomrule
\end{tabular}
\caption*{Table 4.1.1: Details of \textit{Qwen2.5-Math-7B-Instruct}.}
\end{table}
 
\begin{table}[!ht]
\centering
\small
\begin{tabular}{@{}ll@{}}
\toprule
\textbf{Attribute}           & \textbf{Specification}         \\
\midrule
Model name                   & ToRA-7B-v1.0                   \\
Parameter count              & $7\times10^{9}$                \\
Base model                   & LLaMA-2-7B                     \\
Maximum context length       & 4\,096 tokens                  \\
Fine-tuning dataset          & ToRA-Corpus-16k                \\
Release date                 & October 8, 2023                \\
\bottomrule
\end{tabular}
\caption*{Table 4.1.2: Details of \textit{ToRA-7B-v1.0}.}
\end{table}

\paragraph{Prompted Model}  
We used the open-source model \textit{Qwen2.5-Math-7B-Instruct} \cite{yang2024qwen2}, a 7-billion-parameter decoder-only Transformer with an 8\,192-token context window.

\paragraph{Prompt Construction.}  
Thirty algorithmic problems and thirty computational questions (see Table 1) are each paired with nine distinct hinting strategies, yielding \(60\times9=540\) prompts.  The nine strategies are: \emph{Few-shot Hinting}, \emph{Chain-of-Thought (CoT)}, \emph{Automatic CoT (Auto-CoT)}, \emph{Self-Consistency}, \emph{Logical CoT (LogiCoT)}, \emph{Chain-of-Symbol (CoS)}, \emph{Structured CoT (SCoT)}, \emph{ReAct}, and \emph{Clear and Specific Instructions} (for definitions and details, see Appendix~A).

\paragraph{Hint Drafting and Refinement.}  
Initial hints were drafted by prompting \textit{GPT-4} \cite{achiam2023gpt} with: (i) the problem statement, expected answer type, and a true answer (ii) the target hinting strategy, and (iii) an instruction to avoid giving a spoiling hint. The author then reviewed each hint for mathematical correctness and to ensure it guided the model without spoiling the answer (see Appendix A.2 for examples). We note that some of the hints are just a restatement of the corresponding questions in more simplified terms, and the model still achieves high performance on those.

\paragraph{Answer Generation.}  
Each of the 540 prompts was submitted to \textit{Qwen2.5-Math-7B-Instruct}, and the model’s output solutions were collected and saved in the Github repository. 

\paragraph{Answer Evaluation.}  
Generated solutions were compared to the reference true answers using \textit{ToRA-7B-v1.0} \cite{gou2024tora}.  We computed fourteen real-valued scores in \([0,1]\) (Accuracy, Calibration Confidence, Clarity, Coherence, Computational Efficiency, Correctness, Efficiency, Final Accuracy, Pass Rate, Redundancy, Robustness Consistency, Similarity, Stepwise Correctness, Validity Rate), following metric definitions.

\subsection{Part II: Empirical Verification of the Analytic Number Theory Conjecture}

\paragraph{Data Preprocessing.}  
From LMFDB \cite{lmfdb} we retrieve the first 25 nontrivial zeros
\(\rho_j=\tfrac12+i\,\gamma_j\) for each Dirichlet character of modulus \(q\le100\).  The imaginary parts \(\{\gamma_j\}\) are centered and scaled to produce scale-invariant raw input vectors.

\paragraph{Problem Formulation.} 
Concretely, given $k=25$ imaginary parts of the first $k=25$ zeros of a Dirichlet $L$-function, 
\[
(\gamma_1,\dots,\gamma_k)
\;\longmapsto\;
q,
\]
we treat the conjecture as a multiclass classification task over prime moduli \(q\le100\) (see Appendix B for more details).

\paragraph{Feature Engineering.}  
In addition to the raw zero sequence \(\{\gamma_j\}\), we incorporate the following statistical properties of each raw zeros vector in the training process as engineered features:
\noindent
\begin{tabular}{@{}p{2.5cm}@{\quad}p{4.5cm}@{}}
\textbf{Sample moments:} &
  mean, variance, skewness, kurtosis, root-mean-square \\[6pt]
\textbf{Difference statistics:} &
  mean first-difference, variance first-difference, mean second-difference, variance second-difference\\[6pt]
\textbf{FFT features:} &
  magnitudes of the first 30 Fourier components
\end{tabular}

\paragraph{Classification Pipeline.}  
We trained two base learners:
\begin{itemize}
  \item a random forest of 200 trees with balanced class weights;
  \item a LightGBM multiclass classifier with early stopping (1 500 estimators, 127 leaves, early stopping after 75 rounds).
\end{itemize}
Hyperparameters are chosen via cross-validation on the training set.

\section{Experiments $\&$ Results}
\subsection{Part I: LLM Performance on Algorithmic Number-Theory Tasks}

We evaluated \textit{Qwen2.5-Math-7B-Instruct} on thirty classical algorithms and thirty computational questions using nine distinct hinting strategies (see Appendix~A). Across all 540 prompt variants, the model attained uniformly strong scores: mean Accuracy, Similarity, Coherence, Validity Rate, Correctness, and Clarity all lie between 0.80 and 0.91 for algorithmic problems, and between 0.73 and 0.93 for computational questions (see Heatmaps $5.1.1-2$). Remarkably, for every individual problem and question there exists at least one hinting strategy under which the model achieves an Accuracy of 0.95 or higher, illustrating the importance of prompt design (see Appendix~A and the Github repository). Even the lowest-performing metrics remained above 0.70 under a suitable choice of hint, confirming a model's ability to generate both algorithms and computational solutions.

\begin{figure}[!ht]
  \centering
  \includegraphics[width=\linewidth]{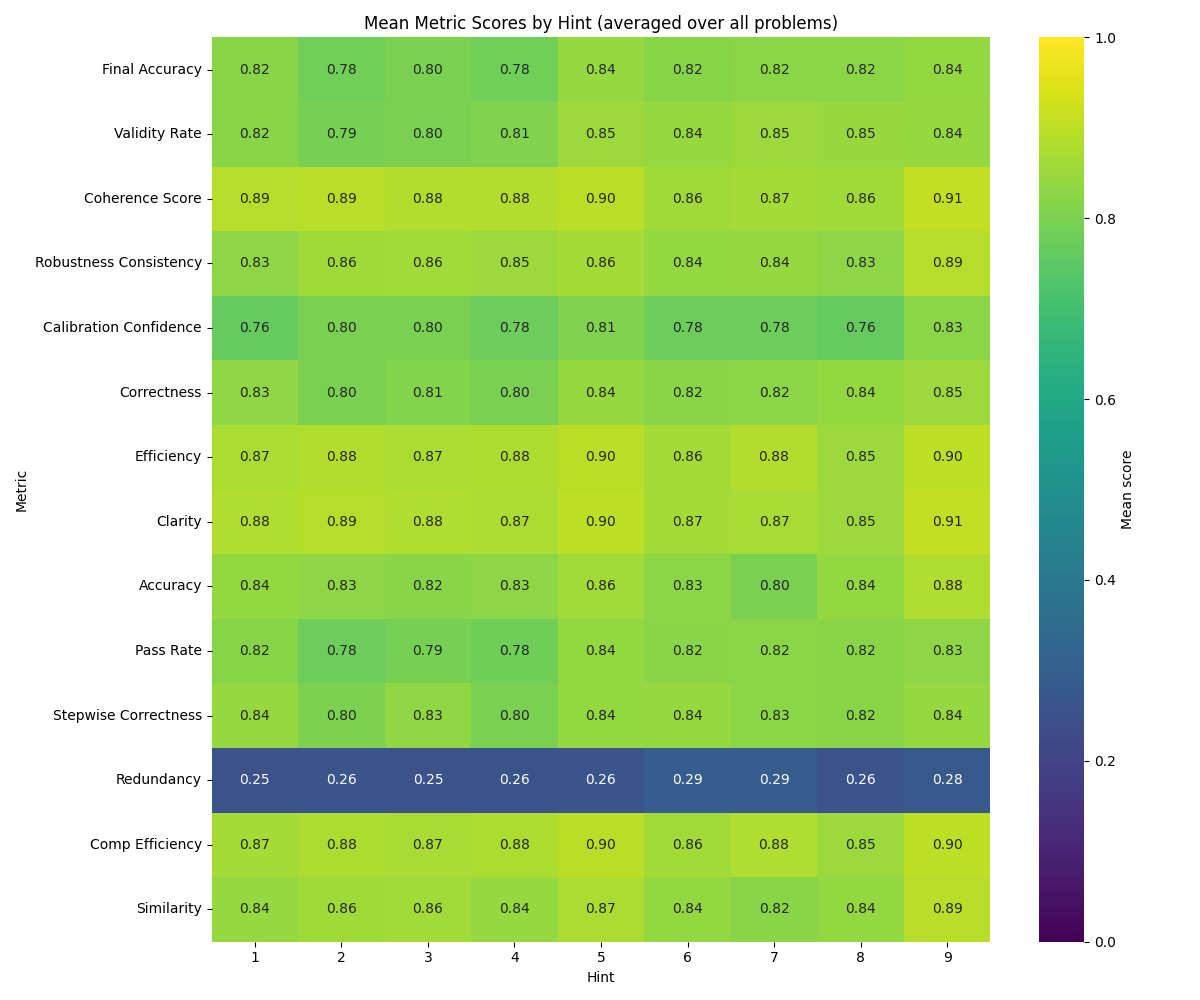}
  \caption*{Heatmap 5.1.1: Mean metric scores by hint across all 30 algorithmic problems (see Appendix A.1).}
\end{figure}

\begin{figure}[!ht]
  \centering
  \includegraphics[width=\linewidth]{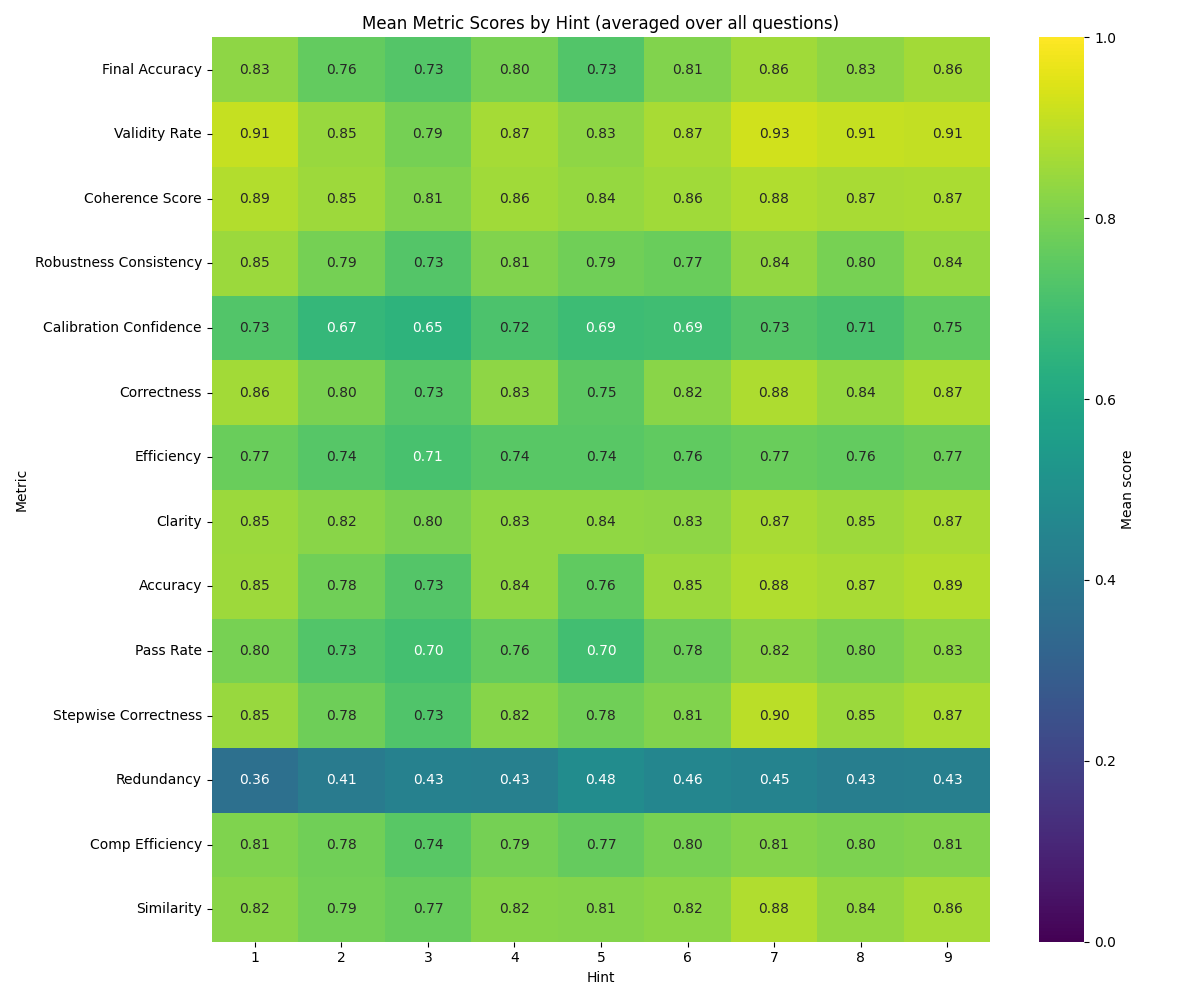}
  \caption*{Heatmap 5.1.2: Mean metric scores by hint across all 30 computational questions (see Appendix A.1).}
\end{figure}
Finally, we controlled output truncation by setting a 1024-token limit for algorithmic prompts and 2048 for computational ones; this introduced modest redundancy mean ( $\approx 0.27$ and $0.44$ respectively), but the redundancy per generated answer was minimal ( $\approx 0.00$) for the right choice of hinting strategy (see Appendix~A).

\subsection{Part II: Empirical Verification of the Analytic Number Theory Conjecture}
As an initial experiment, we trained the random forest on a small dataset of Dirichlet characters—using only the first 25 nontrivial zeros—and evaluated on a held‐out test set of size 21, yielding 61.9\% test accuracy. Augmenting the feature vectors with engineered statistics raised accuracy to 100\%.  Consequently, we expanded to the full range of modulus $5 \le q \le 200$ (214 test points) and incorporated a richer feature set—zero moments, first‐ and second‐difference statistics, and the first 30 FFT magnitudes—into the LightGBM multiclass classifier.  On this larger dataset we achieved 96.7\% validation accuracy and 93.9\% test accuracy.  

\begin{figure}[!ht]
  \centering
  \includegraphics[width=0.5\textwidth]{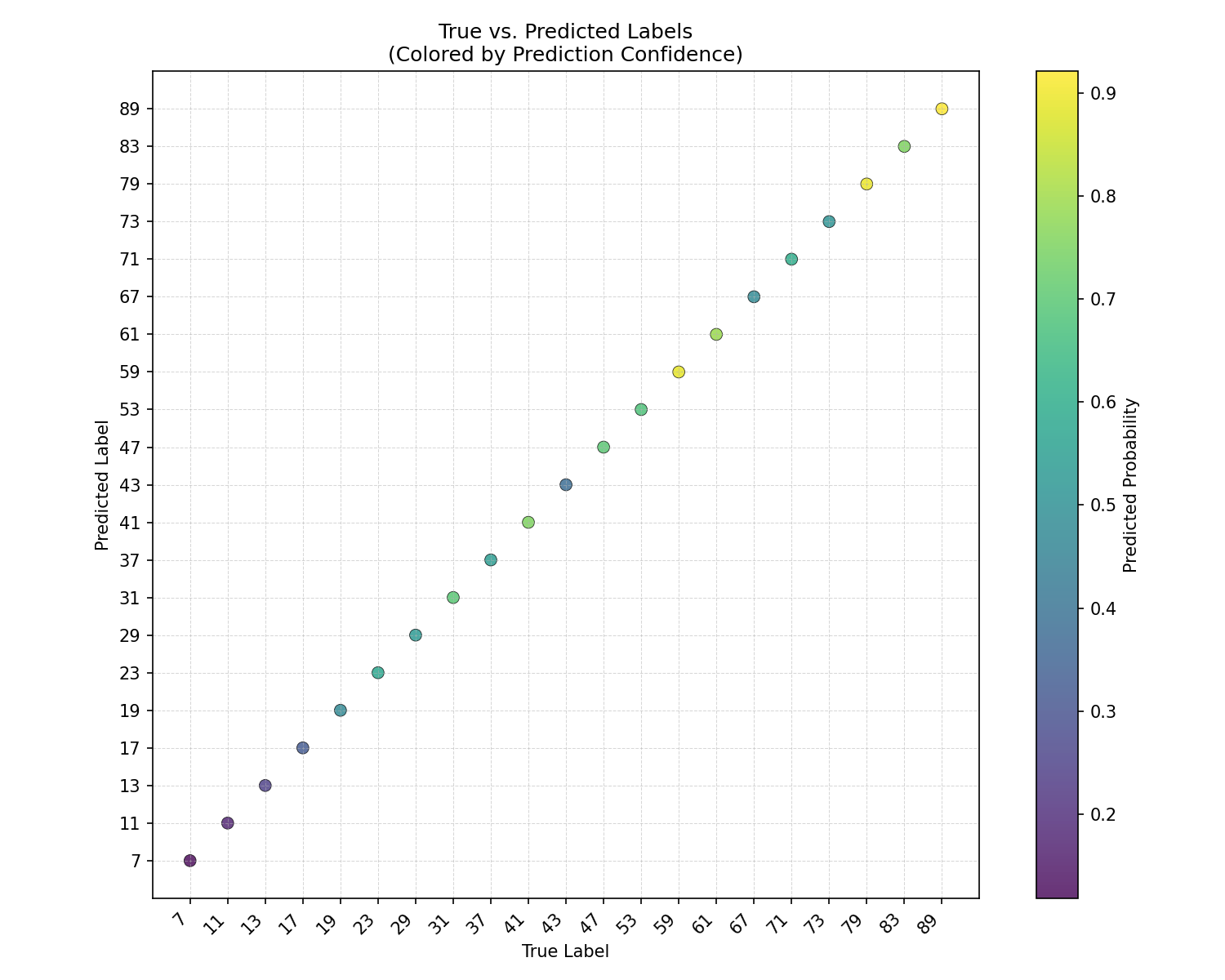}
  \caption*{Heatmap 5.2.1: True vs.\ predicted labels for the engineered‐feature classifier (100\% accuracy). All points lie on the diagonal, indicating perfect classification.}
\end{figure}

\section{Conclusion and Future Directions}

We have presented two complementary applications of LLMs and classical machine learning algorithms to algorithmic and analytic number theory.  First, by constructing the \emph{Hinted Algorithmic Number Theory} (HANT) benchmark of thirty algorithmic problems and thirty computational questions, accompanied by nine non-spoiling hinting strategies, we demonstrated that the state-of-the-art LLM \textit{Qwen2.5‑Math‑7B‑Instruct} can solve every problem/question with at least 95\% accuracy when given an appropriate hint.  Moreover, for a fixed hinting style the model’s mean accuracy reaches 0.88 on algorithms and 0.89 on computations, highlighting the high performance at algorithmic number theory.

In the second study, we treated the folklore conjecture on Dirichlet $L$-function zeros as a multiclass classification problem.  After incorporating carefully engineered statistical features in the training process on labels less than or equal to 200, we could empirically verify the conjecture with a test accuracy of 93.9\%.  The importance of the statistical features in the training process suggests two new conjectures for future studies: that each $L$-function’s zeros encode hidden statistical patterns, and that zeros corresponding to different characters sharing the same modulus possess a similar underlying pattern.

\paragraph{Future Directions.}  
A natural next step is to fine‐tune a mathematics‐focused language model on datasets similar to the HANT dataset. On the analytic side, there is a natural follow-up: one can attempt to provide more empirical or theoretical justification for the two proposed conjectures.
\section*{Acknowledgments}
This work was carried out without any external funding. The author would like to thank Ghaith Hiary (Ohio State University) and Zhihui Zhu (Ohio State University) for answering a few brief questions.

\section*{Appendix}
The appendix comprises two main sections (A and B). In Appendix A, we establish our experimental setup for evaluating a state-of-the-art LLM on algorithmic number‐theory tasks. Section A.1 details the experiments and the evaluations of model's outputs, highlighting a strong performance, while section A.2 presents a catalogue of problem-hint (respectively, question-hint) pairs for a selection of the thirty problems (respectively, thirty question types) listed in Table 1, with model’s generated solutions for observation and analysis. Appendix B is devoted to the second part of the paper on the empirical verification of the analytic number theory conjecture. It begins by presenting the background, the precise statement of the conjecture, and the experimental design. In Section B.1, we show that training a random forest solely on the raw sequence of zeros yields a relatively low test accuracy. Since the conjecture is widely believed true, this result suggests the existence of additional statistical structure in the zeros beyond their consecutive spacings. To address this, we demonstrate that incorporating carefully chosen statistical properties of the zeros (mean spacing, variance, skewness, etc.) into the training process of the random forest leads to perfect ($100 \%$) test accuracy. We note that studies in section B.1 are restricted to only 21 labels. Lastly, section B.2 shows that a LightGBM classifier achieves 0.939 test accuracy on a test dataset of size 214 (140 different labels) when sufficient statistics of the zeros are incorporated in the training process.  
 \begin{table}[!ht]
  \small                            
  \setlength{\tabcolsep}{4pt}       
  \begin{tabularx}{\columnwidth}{@{}>{\raggedright\arraybackslash}p{0.30\columnwidth}
                                  >{\raggedright\arraybackslash}X@{}}
    \toprule
    \textbf{Category} & \textbf{Algorithms/ Types of Questions}\\
    \midrule
    GCD \& Linear Congruence &
      Euclid’s Algorithm (2); Extended Euclidean Algorithm (13); Binary GCD (23);\\
    & Multiplicative Inverse (3); Linear Congruence Solver (16); Chinese Remainder Theorem (15)\\
    \addlinespace
    Diophantine \& Continued Fractions &
      Solving Pell’s Equation (1); Zeckendorf’s Representation (22);\\
    & Continued Fraction Expansion (4)\\
    \addlinespace
    Integer Representation &
     Greedy Egyptian Fractions Algorithm (18); Rational Root Theorem Algorithm (17); Fermat’s Two-Squares Theorem (24)\\
    \addlinespace
    Modular-Arithmetic \& Root-Extraction &
      Hensel-Lifting of Cube Roots (5); Primitive Roots mod \(p\) (6); Square Roots mod \(p\) (7)\\
    \addlinespace
    Exponentiation \& Discrete Log &
      Binary Exponentiation (8); Discrete Logarithm mod \(p\) (25)\\
    \addlinespace
    Primality Testing &
      Primality Testing via Trial Division (9); Primality Testing Using Euler’s Criterion (10);\\
    & Fermat’s Primality Test (20); Lucas–Lehmer Primality Test (29)\\
    \addlinespace
    Sieves \& Totient-Functions &
      Sieve of Eratosthenes (11); Euler’s Totient Sieve (21); Carmichael Function Computation (19)\\
    \addlinespace
    Factorization \& Ideals &
      Fermat Factorization Method (14); Decomposition of a Prime Ideal (30)\\
    \addlinespace
    Cryptographic-Key \& EC Operations &
      RSA Key Generation (27); Elliptic Curve Point Multiplication (26); Elliptic Curve Point Doubling (28)\\
    \bottomrule
  \end{tabularx}
  \caption*{Table 1: Famous number-theoretic algorithms and computational questions drawn from classical textbooks on algorithmic and computational number theory \cite{cohen2013course,shoup2009computational}. Each problem type is followed by its index, which we follow in this paper.}
\end{table}

\subsection*{Appendix A: On the Performance of an open-source state-of-the-art LLM on Algorithmic Number Theory}

The model \textit{Qwen2.5‑Math‑7B‑Instruct} \cite{yang2024qwen2} was prompted on 30 algorithmic problems and 30 corresponding computational questions in algorithmic number theory for algorithms and computational answers generation, respectively. For each problem and each question, 9 prompting strategies (more precisely, hinting strategies) were employed, whence the total number of prompts of the model is $60 \times 9= 540$ prompts. Each prompt fixed the problem or question and varied only the hint according to one of the following strategies: 
\begin{enumerate} 
    \item  \textbf{Few-shot Hinting}: A few examples were provided to the model as a hint. 
    \item \textbf{Chain-of-thought (CoT) Prompting}: This presents a hint about the intermediate steps in the derivation process of the answer.
    \item \textbf{Automatic Chain-of-Thought (Auto-CoT) Prompting}: This presents the hint in a step-by-step way of thinking.
    \item \textbf{Self-Consistency}: This shows some different reasoning paths about the problem or question and selects the most consistent one. 
    \item \textbf{Logical Chain-of-Thought (LogiCoT) Prompting}: This presents the hint in a sequence of reasoning and logical verification steps.
    \item \textbf{Chain-of-Symbol (CoS) Prompting}: This presents the hint mathematically and using mathematical symbols.
    \item \textbf{Structured Chain-of-Thought (SCoT) Prompting}: This presents the hint for the problem or question in a programming way, e.g., the use of loops, equations on different lines, etc. $\cdot$ 
    \item \textbf{ReAct Prompting}: The hint is presented in a sequence of action commands and thought steps.
    \item \textbf{Clear and Specific Instructions}: The hint is presented in a sequence of precise and clear instructions.
\end{enumerate}
The hints were initially drafted by prompting \textit{GPT-4} \cite{achiam2023gpt}  with three inputs: (i) the problem statement (including the expected answer type—algorithm or computational solution), (ii) the designated hinting strategy, and (iii) an instruction to avoid revealing the full solution. Subsequently, the author reviewed and refined the hints to ensure that they are mathematically correct and do not spoil the answers (see \textit{Appendix A.2} for examples and further details). 

Generated answers for each problem–hint (or question-hint) pair were evaluated by the \textit{ToRA-7B-v1.0} model \cite{gou2024tora} using fourteen real‐valued performance metrics in $[0,1]$: Accuracy, Calibration Confidence, Clarity, Coherence Score, Computational Efficiency, Correctness, Efficiency, Final Accuracy, Pass Rate, Redundancy, Robustness Consistency, Similarity, Stepwise Correctness, and Validity Rate. It is important to mention that the true answers (up to reformulation and little changes) corresponding to each problem and question statement were given as part of the evaluation prompts to the model \textit{ToRA-7B-v1.0} for comparison and evaluation. We define the Similarity metric as the percentage of similarity between a generated answer and the corresponding actual "true" answer. Other metric definitions follow.

\subsubsection*{Appendix A.1: Performance of the LLM on the Generation of Algorithms}
Heatmaps $A.1.1$ and $A.1.2$ show the means of the metric scores averaged over all thirty algorithmic problems and all thirty computational questions, respectively. Remarkably, the model achieves very high performance on both task types: the mean scores for Accuracy, Similarity, Coherence, Validity Rate, Correctness, and Clarity lie in the range \([0.80,0.91]\) for algorithmic problems and in the range \([0.73,0.93]\) for computational questions. We observe that accuracy for algorithmic problems peaks at $0.88$ with Hint $9$ (Clear and Specific Instructions), while the maximum Validity Rate ($0.93$) was achieved by Hint 7 (Structured Chain-of-Thought (SCoT) Prompting). In fact, apart from the Redundancy metric, all other metrics remained in similar ranges, with the worst numbers fluctuating around $0.7$. Speaking of Redundancy, we note that the model was prompted in a way that the generated outputs contain the problem and hint statements in addition to the generated solution.  The author tested prompting the model with different allowed "number of tokens", and it turned out, with a small that, with a modest "maximum number of tokens," the generated solutions were truncated midway for some of the algorithmic problems and computational questions (in some cases, computational solutions required a higher number of tokens). Thus, the author decided to set 
\begin{lstlisting}
  "algorithmic_max_tokens": 1024,
  "computational_max_tokens": 2048
\end{lstlisting}
Thus, achieving strong performance came at the cost of some redundancy (around $0.27$ for algorithmic solutions and around $0.44$ for computational solutions), i.e., some generated outputs were more verbose or repetitive, but the majority was not.
\begin{figure}[!ht]
    \centering
    \includegraphics[width=\linewidth]{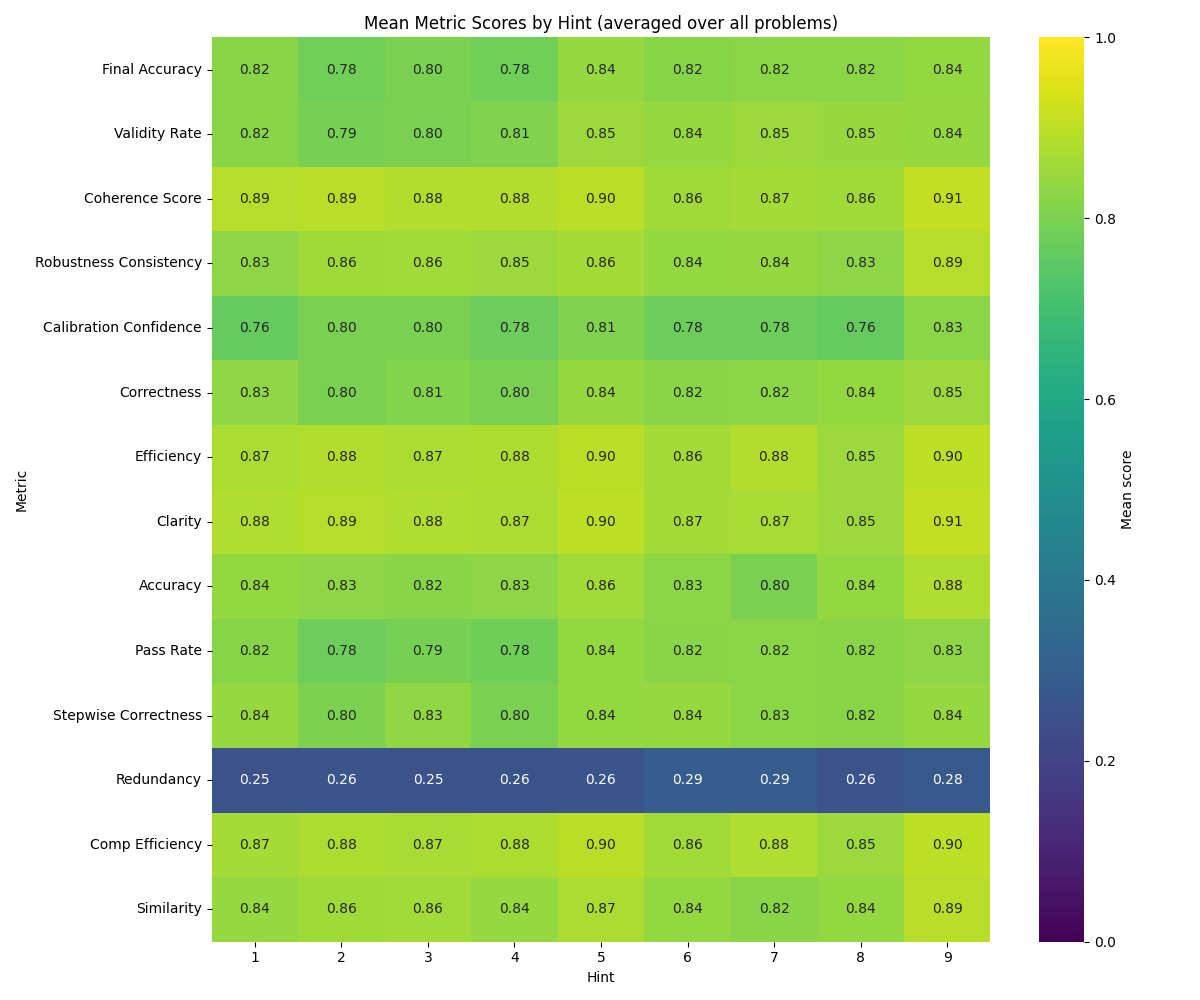}
    \caption*{Heatmap A.1.1: \textbf{Mean Metric Scores by Hint}: The heatmap presents the means of the metric scores averaged over all 30 \textbf{algorithmic problems}.}
\end{figure}
\begin{figure}[!ht]
    \centering
    \includegraphics[width=\linewidth]{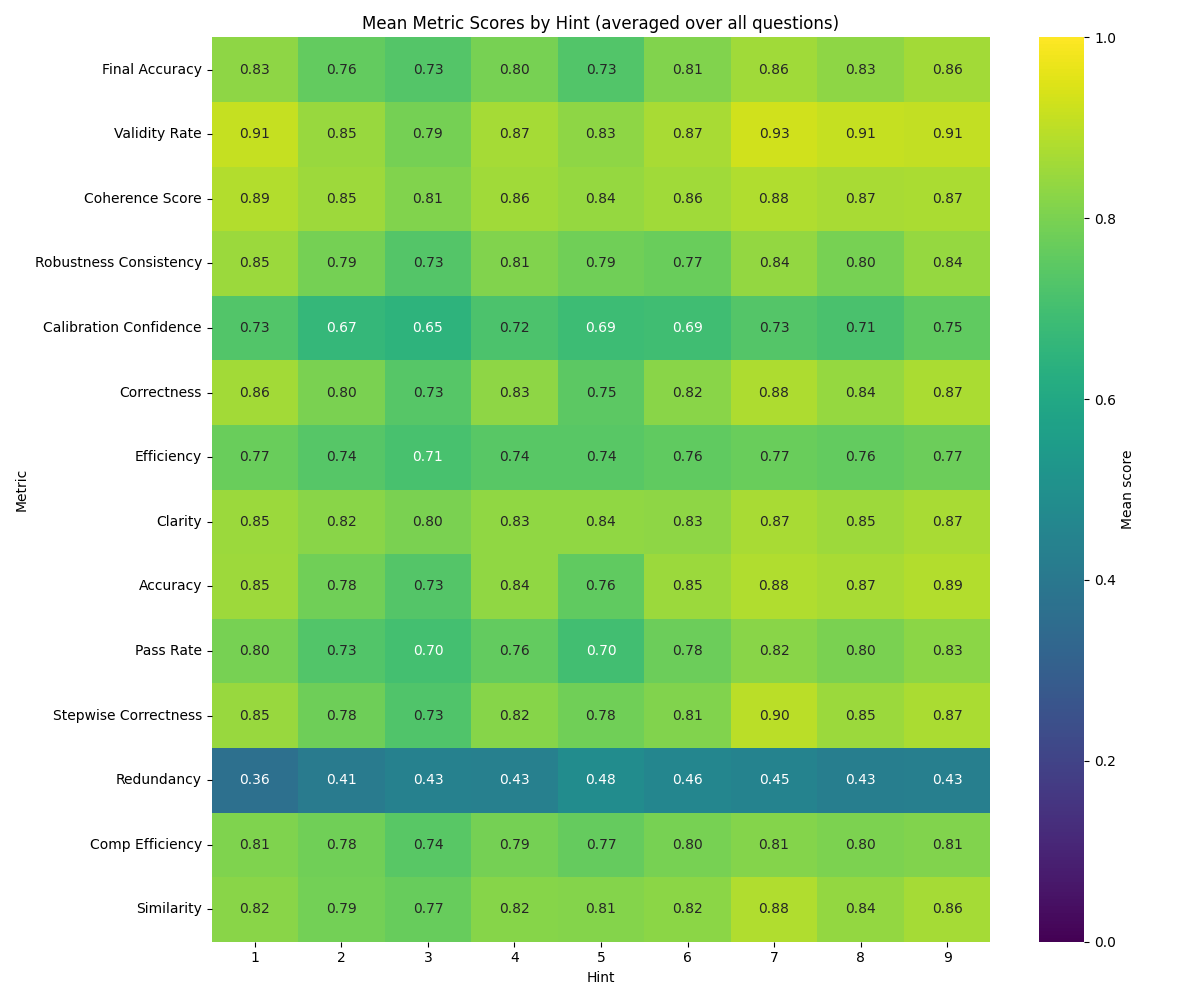}
    \caption*{Heatmap A.1.2 \textbf{Mean Metric Scores by Hint}: The heatmap presents the means of the metric scores averaged over all 30 \textbf{computational questions}.}
\end{figure}

Figures $A.1.3-A.1.8$ present 3D surfaces of the Accuracy, Correctness, and Similarity metrics, where the vertical axis encodes the score in $[0,1]$ over the $30$ problems or questions ($y$-axis) and nine hints ($x$-axis). Figure $A.1.3$ reveals that the LLM achieves near-perfect accuracy scores (score $\ge 0.9$) on the majority of algorithmic problem-hint pairs (on approximately $75 \%$ of the pairs), while a modest accuracy ($[0.6, 0.8]$) for all problems under a suitable choice of a hinting strategy (hints $2-5$). However, some problem-hint pairs have poor accuracy score, e.g, for problems $21-24$, some problem-hint pairs have accuracy score $\le 0.4.$ Similarly, as illustrated in Figures $A.1.4$ and $A.1.5$, the correctness and similarity scores on algorithmic problems exhibit high scores $[0.75, 1]$ for almost all hinting strategies for most of the problems, while some problems attain low scores (score $\le 0.5$) for some hinting strategies, but all problems could be solved with correctness and similarity scores at least $0.7$ for an appropriate choice of hint. On the other hand, the accuracy surface for computational questions, Figures $A.1.6$, forms a broad plateau (\(\approx0.65\!-\!0.85\)) with peaks (\(\ge0.9\)) achieved by question-hint pairs for hints 3–7 align and questions 6–18, and troughs (\(\le0.45\)) for some other pairs. For instance, for questions $23-26$, accuracy score stabilized at 0 for all hinting strategies except for a few hinting strategies, namely hint $5-9 $, where accuracy score rise to near $1.0$ (this is clear by the gap in the surface). However, there is one exception, question $26,$ where accuracy score never exceed $0.5$. These observations can be confirmed using the evaluation tables for each question published on Github. The correctness surface, Figures $A.1.7$, is essentially the same shape as the accuracy surface. However, the similarity score, Figures $A.1.8$, is mostly elevated (\(\approx0.80\!-\!1.00\)) and in the worst case falls in (\(\approx0.60\!-\!0.75\)) across almost all hint–question pairs. 

\newpage
\clearpage
\begin{figure}[!ht]
    \centering
    \includegraphics[width=\linewidth]{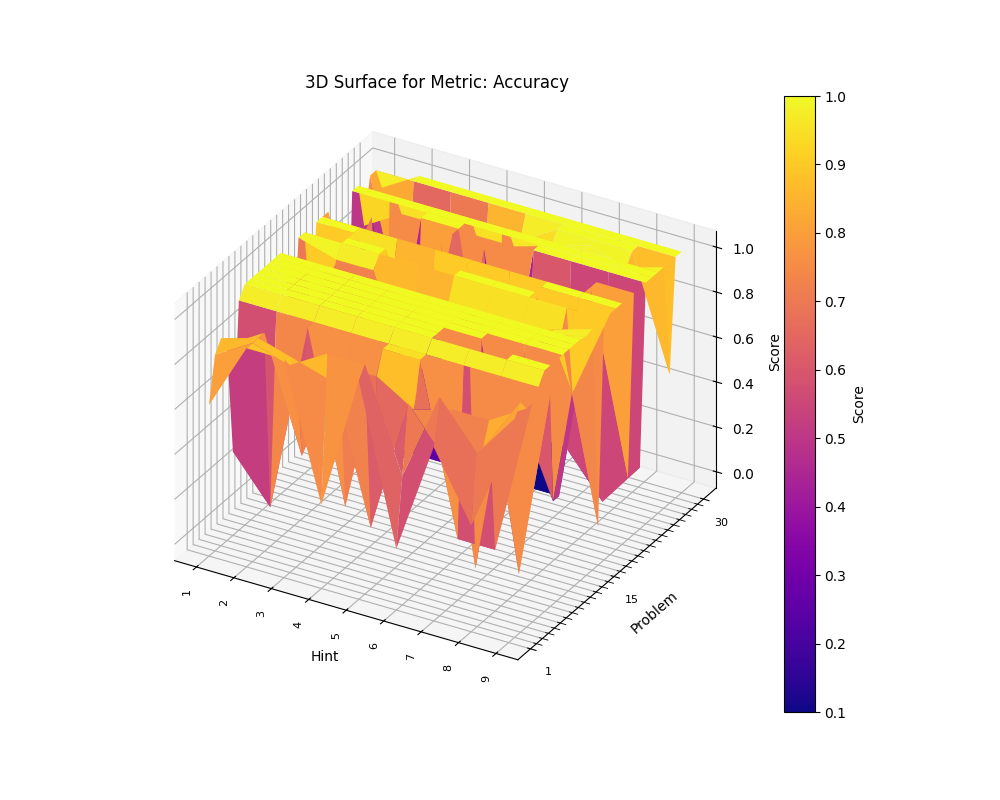}
    \caption*{Figure A.1.3: \textbf{3D surface of the Accuracy metric for algorithmic problems} }
\end{figure}

\begin{figure}[!ht]
    \centering
    \includegraphics[width=\linewidth]{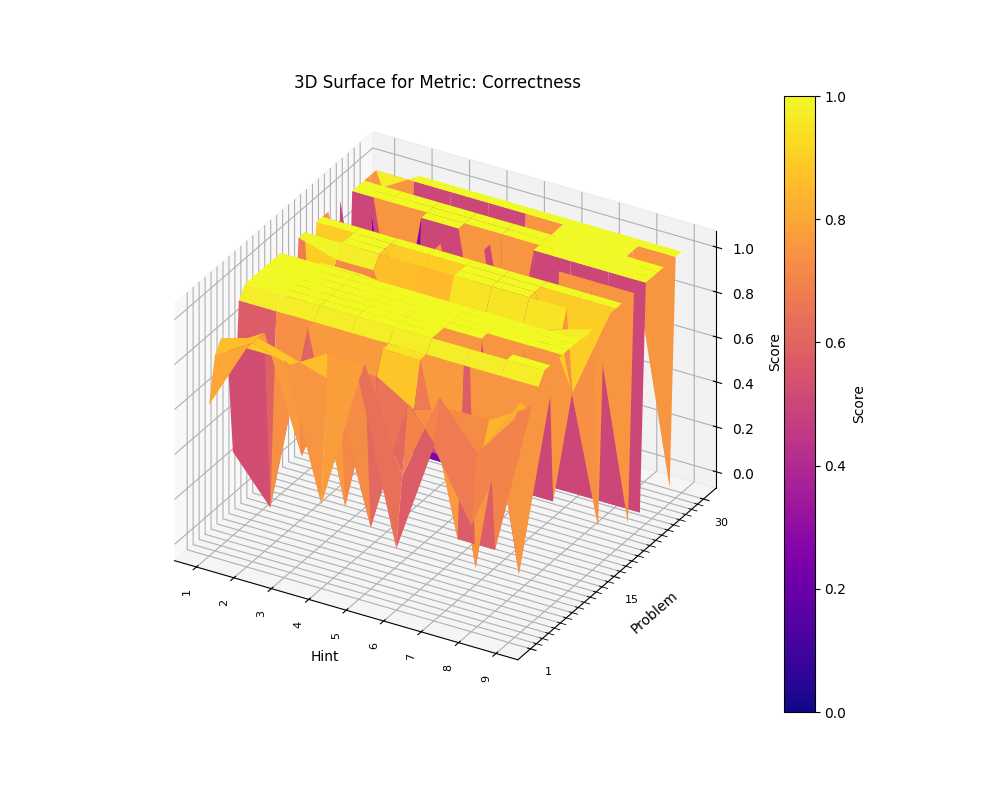}
    \caption*{Figure A.1.4: \textbf{3D surface of the Correctness metric for algorithmic problems}}
\end{figure}

\begin{figure}[!ht]
    \centering
    \includegraphics[width=\linewidth]{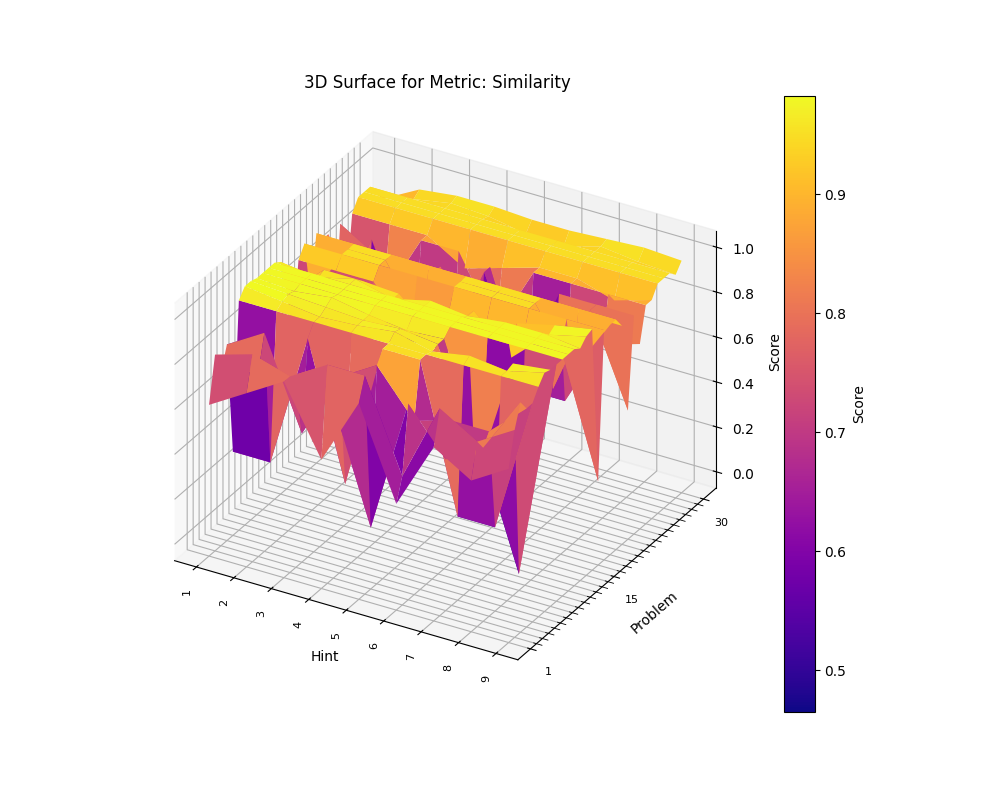}
    \caption*{Figure A.1.5: \textbf{3D surface of the Similarity metric for algorithmic problems}}
\end{figure}

\begin{figure}[!ht]
    \centering
    \includegraphics[width=\linewidth]{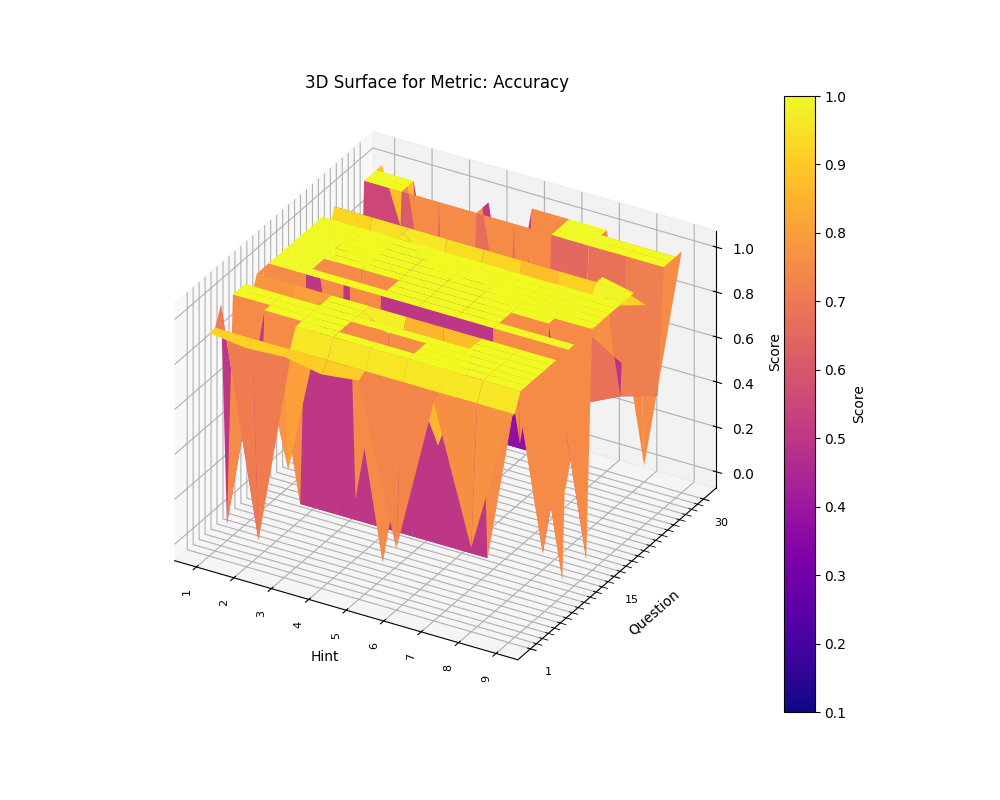}
    \caption*{Figure A.1.6: \textbf{3D surface of the Accuracy metric for computational questions}}
\end{figure}

\begin{figure}[!ht]
    \centering
    \includegraphics[width=\linewidth]{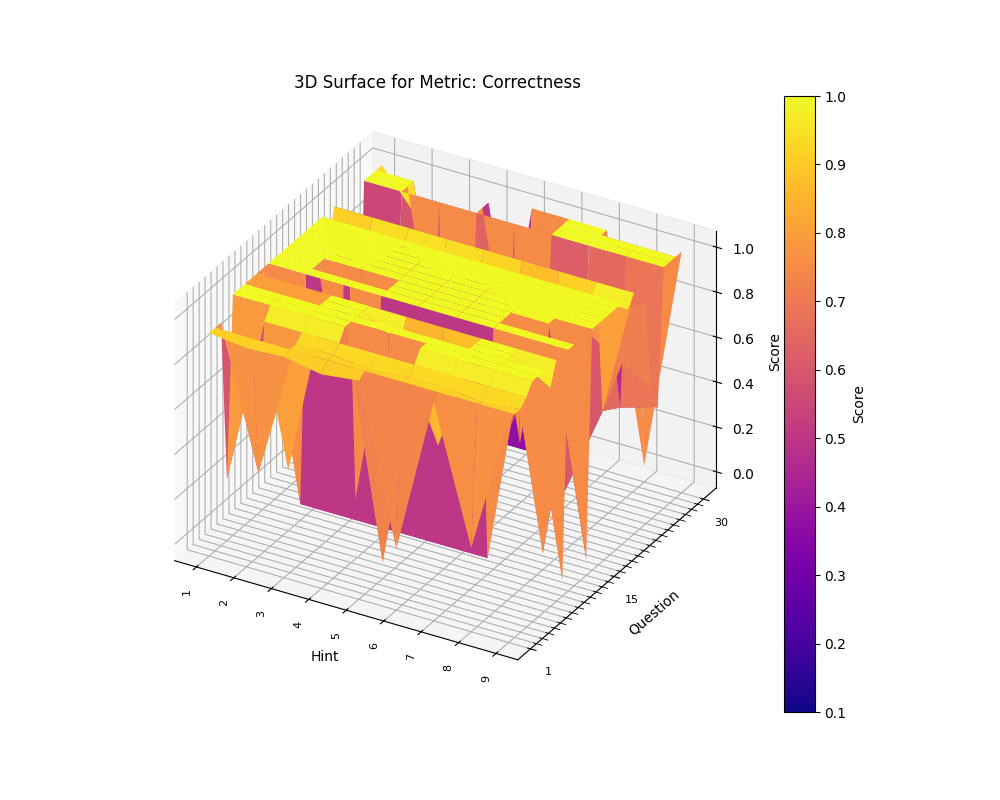}
    \caption*{Figure A.1.7: \textbf{3D surface of the Correctness metric for computational questions}}
\end{figure}

\begin{figure}[!ht]
    \centering
    \includegraphics[width=\linewidth]{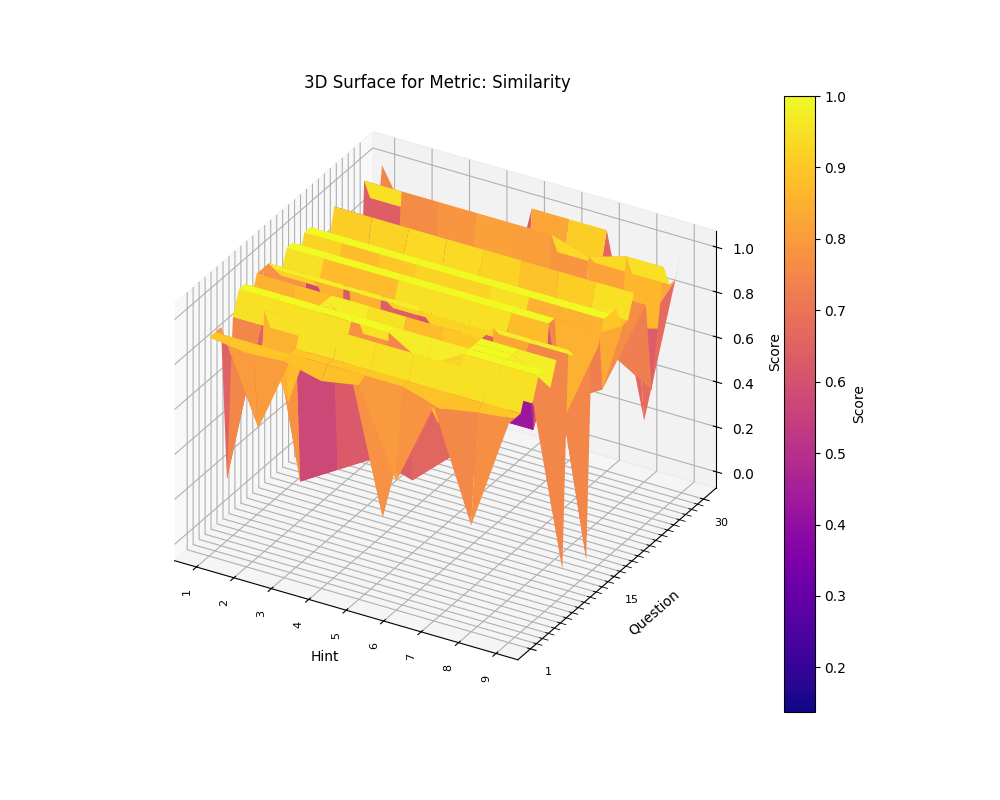}
    \caption*{Figure A.1.8: \textbf{3D surface of the Similarity metric for computational questions}}
\end{figure}
\clearpage
\newpage
\subsubsection*{Appendix A.2: Highly Accurate Generated Outputs for Algorithmic Problems and Computational Questions}

To illustrate the model’s overall behavior, this section presents a representative sample of generated outputs for selected problem–hint and question–hint pairs from each category in Table 1.  For readers interested in a more detailed view, the complete set of outputs is available in the linked GitHub repository. Given the high performance of the model, there are only a very few pairs of problem-hint and question-hint that have partially accurate solutions, a sample of which we present in section \textit{Appendix A.3}. In this section, however, we present the generated solutions that are highly accurate. For each algorithmic problem and computational question below, we report only the Accuracy metric; readers can consult the GitHub repository for all other metrics, which exhibit behavior highly similar to the Accuracy metric. Also, the actual correct solution to an algorithmic problem or a computational question, $n$, is available in the subsection named "Algoritm $/$ Answer" in the corresponding file named "Problem $n$" or "Question $n$", respectively.

\subsubsection*{Solving Pell’s Equation (Computational Question 1)- Self-Consistency; Accuracy: $0.95$}
\fbox{%
  \begin{minipage}{\dimexpr\linewidth-2\fboxsep-2\fboxrule\relax}
    \scriptsize
    \noindent\textbf{Question.}
    Find the fundamental solution \((x,y)\) to Pell’s equation
    \[
      x^2 - 3\,y^2 = 1,
    \]
    for integers \((x,y)\neq(1,0)\).

    \medskip
    \noindent\textbf{Hint 4 (Self-Consistency).}
    \begin{itemize}[leftmargin=1.5em, topsep=0pt, itemsep=0pt]
      \item Make sure you understand what continued fractions are.
      \item Note \(x^2 - d\,y^2=1\) implies 
        \(\tfrac{x}{y}=\sqrt{d+\tfrac1{y^2}}\approx\sqrt{d}\) for large \(y\).
      \item Compare your result with known small examples.
    \end{itemize}

    \medskip
    \noindent\textbf{Solution.}
    The continued fraction of \(\sqrt{3}\) is
    \[
      \sqrt{3} = [1;\overline{1,2}],
    \]
    whose convergents are
    \(\tfrac11,\;\tfrac21,\;\tfrac53,\;\dots\).  Checking
    \[
      2^2 - 3\cdot1^2 = 4 - 3 = 1
    \]
    shows that the fundamental solution is
    \[
      \boxed{(2,\,1)}.
    \]
  \end{minipage}%
}

\bigskip

\subsubsection*{Solving Pell’s Equation (Algorithmic Problem 1)- Logical Chain-of-Thought Prompting; Accuracy: $0.9$}
\fbox{%
  \begin{minipage}{\dimexpr\linewidth-2\fboxsep-2\fboxrule\relax}
    \scriptsize
    \noindent\textbf{Problem Statement.}
    Given a nonsquare positive integer \(d\), solve
    \[
      x^2 - d\,y^2 = 1
    \]
    for integers \((x,y)\neq(1,0)\).

    \medskip
    \noindent\textbf{Hint 5 (Logical CoT).}
    \begin{itemize}[leftmargin=1.5em, topsep=0pt, itemsep=0pt]
      \item The relation \(x^2 - d\,y^2 = 1\) means \(x/y\) is nearly \(\sqrt{d+\frac{1}{y^2}}\) and for large \(y\), close to \(\sqrt{d}\).
      \item Recall continued fractions.
      \item Verify that your solution satisfies the original equation.
    \end{itemize}

    \medskip
\noindent\textbf{Generated Answer.}

To solve Pell's equation \(x^2 - d\,y^2 = 1\) using the continued fraction expansion of \(\sqrt{d}\), we follow these steps:

\begin{enumerate}[leftmargin=2em, topsep=0pt, itemsep=0pt]
  \item \textbf{Compute the continued fraction expansion of \(\sqrt{d}\):}\\
    The continued fraction expansion of \(\sqrt{d}\) is given by:
    \[
      \sqrt{d} = [a_0; \overline{a_1, a_2, \ldots, a_{2m}}]
    \]
    where \(a_0 = \lfloor \sqrt{d} \rfloor\) and the sequence \(a_1, a_2, \ldots, a_{2m}\) is periodic with period \(2m\).
  \item \textbf{Generate the convergents of the continued fraction:}\\
    The convergents of the continued fraction are given by:
    \[
      \frac{p_n}{q_n} = [a_0; a_1, a_2, \ldots, a_n]
    \]
    for \(n = 0, 1, 2, \ldots, 2m\). They satisfy the recurrences:
    \[
      p_n = a_n p_{n-1} + p_{n-2}, 
      \quad
      q_n = a_n q_{n-1} + q_{n-2},
    \]
    with \(p_{-2} = 0,\ p_{-1} = 1,\ q_{-2} = 1,\ q_{-1} = 0\).
  \item \textbf{Identify the fundamental solution:}\\
    The fundamental solution \((x, y)\) to \(x^2 - d\,y^2 = 1\) is given by the convergent \(\tfrac{p_{2m-1}}{q_{2m-1}}\) if the period length is even.
  \item \textbf{Verify the solution:}\\
    Check that
    \[
      p_{2m-1}^2 - d\,q_{2m-1}^2 = 1.
    \]
\end{enumerate}

\[
\boxed{
\begin{aligned}
&1.\ \text{Compute the continued fraction expansion of }\sqrt{d}\text{ as } \\ & [a_0;\overline{a_1,\ldots,a_{2m}}].\\
&2.\ \text{Generate convergents }p_n/q_n\text{ via }p_n=a_np_{n-1}+p_{n-2}, \\ & \ q_n=a_nq_{n-1}+q_{n-2},\ p_{-2}=0,p_{-1}=1,q_{-2}=1,q_{-1}=0.\\
&3.\ \text{Fundamental solution: }(x,y)=(p_{2m-1},q_{2m-1}).\\
&4.\ \text{Verify }p_{2m-1}^2 - d\,q_{2m-1}^2 = 1.
\end{aligned}
}
\]

  \end{minipage}%
}

\subsubsection*{Multiplicative Inverse (Computational Question 3) — ReAct Prompting; Accuracy: \(1.0\)}
\fbox{%
  \begin{minipage}{\dimexpr\linewidth-2\fboxsep-2\fboxrule\relax}
    \scriptsize
    \noindent\textbf{Problem Statement.}
    Given the integers \( a = 3 \) and \( n = 11 \), compute the multiplicative inverse of 3 in \(\mathbb{Z}/11\mathbb{Z}\).
    In other words, find an integer \( x \) such that:
    \[
      3 \cdot x \equiv 1 \pmod{11}.
    \]
    Provide all the intermediate steps.

    \medskip
    \noindent\textbf{Hint 8 (ReAct Prompting).} \\
    1. **Action**: Use the extended Euclidean algorithm to find an equation between the GCD, \( a \) and \( n \).   \\
    2. **Reasoning**: If the GCD is 1, this seems to be close to the definition of inverse.   \\
    3. **Reasoning**: If the GCD is not 1, what happens—can you say anything about the left and right hand side of the equations? \\
    4. **Action**: Provide the computational answer \\

    \medskip
    \noindent\textbf{Generated answer.}

    To find the multiplicative inverse of \( a = 3 \) in \( \mathbb{Z}/11\mathbb{Z} \), we need to find an integer \( x \) such that:
    \[
      3x \equiv 1 \pmod{11}.
    \]

    We will use the Extended Euclidean Algorithm to solve this. The Euclidean Algorithm is used to find the greatest common divisor (GCD) of two integers, and the Extended Euclidean Algorithm also provides the coefficients (in this case, \( x \)) of Bézout’s identity, which states:
    \[
      \gcd(a, n) = a\,x + n\,y.
    \]

    Apply the Euclidean Algorithm to \( a = 3 \) and \( n = 11 \):
    \begin{enumerate}[leftmargin=2em, topsep=0pt, itemsep=0pt]
      \item \(11 = 3 \cdot 3 + 2\).
      \item \(3  = 2 \cdot 1 + 1\).
      \item \(2  = 1 \cdot 2 + 0\).
    \end{enumerate}
    Since the remainder is 0, \(\gcd(3,11)=1\), so an inverse exists.

    Work backwards to express 1 as a combination of 3 and 11:
    \[
      1 = 3 - 2\cdot1,
      \quad
      2 = 11 - 3\cdot3
      \;\Longrightarrow\;
      1 = 3 - (11 - 3\cdot3) = 3\cdot4 - 11.
    \]
    Thus
    \[
      1 = 3\cdot4 + (-1)\cdot11
      \;\Longrightarrow\;
      3\cdot4 \equiv 1 \pmod{11}.
    \]
    Therefore, the multiplicative inverse of 3 in \(\mathbb{Z}/11\mathbb{Z}\) is
    \[
      \boxed{4}.
    \]
  \end{minipage}%
}

\subsubsection*{Multiplicative Inverse (Algorithmic Problem 3) — Chain-of-Thought Prompting; Accuracy: \(0.90\)}
\fbox{%
  \begin{minipage}{\dimexpr\linewidth-2\fboxsep-2\fboxrule\relax}
    \scriptsize
    \noindent\textbf{Problem Statement.}
    Given integers \( a \) and \( n \) such that \(\gcd(a,n)=1\), find the multiplicative inverse of \( a \) in \(\mathbb{Z}/n\mathbb{Z}\)—that is, find an integer \( x \) satisfying:
    \[
      a \cdot x \equiv 1 \pmod{n}.
    \]
    \noindent\textbf{Hint 2 ((CoT) Prompting)} 
    \begin{itemize}
        \item Apply the extended Euclidean algorithm to \( a \) and \( n \) to find the GCD.
        \item Can the GCD be bigger than 1? If it's 1, what equation do we have?
        \item Deduce the multiplicative inverse modulo \( n \).
    \end{itemize}
   \noindent\textbf{Generated Answer.}

To find the multiplicative inverse of \(a\) in \(\mathbb{Z}/n\mathbb{Z}\), we use the extended Euclidean algorithm. The algorithm not only finds the greatest common divisor (gcd) of \(a\) and \(n\), but also expresses this gcd as a linear combination of \(a\) and \(n\). Since \(\gcd(a, n) = 1\), we can express 1 as a linear combination of \(a\) and \(n\), which gives us the multiplicative inverse of \(a\) modulo \(n\).

Here is the step-by-step algorithm:

\begin{enumerate}[leftmargin=2em, topsep=0pt, itemsep=0pt]
  \item \textbf{Apply the Euclidean Algorithm to find \(\gcd(a,n)\):}
    \[
      n = q_0 a + r_0,\quad
      a = q_1 r_0 + r_1,\]
    \[
      r_0 = q_2 r_1 + r_2,\ \ldots,\ 
      r_{k-2} = q_k r_{k-1} + r_k,\quad
      r_{k-1} = q_{k+1} r_k + 0,
    \]
    where \(r_k = \gcd(a,n)=1\).
  \item \textbf{Express 1 as a linear combination of \(a\) and \(n\):}\\
    Work backwards from
    \[
      r_k = r_{k-2} - q_k r_{k-1},
      \quad
      r_{k-1} = r_{k-3} - q_{k-1} r_{k-2},
    \]
    to obtain
    \[
      r_k = (1 + q_k q_{k-1})\,r_{k-2} - q_k\,r_{k-3},
    \]
    and continue until
    \[
      1 = x\,a + y\,n.
    \]
  \item \textbf{Identify the inverse:}\\
    The coefficient \(x\) of \(a\) in \(1 = x\,a + y\,n\) satisfies
    \[
      x\,a \equiv 1 \pmod{n},
    \]
    so the multiplicative inverse of \(a\) modulo \(n\) is
    \[
      \boxed{x}.
    \]
\end{enumerate}

  \end{minipage}%
}

\subsubsection*{Fermat’s Two-Squares Theorem (Computational Question 24) — Structured Chain-of-Thought (SCoT) Prompting; Accuracy: 0.9}
\fbox{%
  \begin{minipage}{\dimexpr\linewidth-2\fboxsep-2\fboxrule\relax}
    \scriptsize
    \noindent\textbf{Problem Statement.}
    Determine whether the number \( N = 130 \) can be written in the form \( x^2 + y^2 \) with \( x, y \in \mathbb{Z} \). \\ \\
    \noindent\textbf{Hint 7 ((SCoT) Prompting).}
    \begin{enumerate}
      \item Factorize \( N \) into its prime factors.
      \item Distinguish between the two types of odd primes and two types of exponents.
      \item Conclude whether such a representation exists, and optionally indicate how one might compute it.
    \end{enumerate}

    \medskip
    \noindent\textbf{Generated answer.}
    To determine whether the number \( N = 130 \) can be written in the form \( x^2 + y^2 \) with \( x, y \in \mathbb{Z} \), we will use Fermat's Two-Squares Theorem. This theorem states that a positive integer \( N \) can be expressed as the sum of two squares if and only if every prime factor of \( N \) of the form \( 4k+3 \) appears with an even exponent in the prime factorization of \( N \).

    **Step 1: Factorize \( N \) into its prime factors.**

    \[
    N = 130
    \]
    \[
    130 = 2 \times 65
    \]
    \[
    65 = 5 \times 13
    \]
    \[
    130 = 2 \times 5 \times 13
    \]

    **Step 2: Distinguish between the two types of odd primes based on their residue modulo 4.**

    - \( 2 \) is not an odd prime.
    - \( 5 \equiv 1 \pmod{4} \)
    - \( 13 \equiv 1 \pmod{4} \)

    **Step 3: Assess the exponents of the primes in the critical residue class.**

    - The prime \( 2 \) appears with an exponent of 1.
    - The prime \( 5 \) appears with an exponent of 1.
    - The prime \( 13 \) appears with an exponent of 1.

    **Step 4: Determine if the configuration permits a representation as \( x^2 + y^2 \).**

    According to Fermat's Two-Squares Theorem, for \( N \) to be expressible as \( x^2 + y^2 \), every prime factor of the form \( 4k+3 \) must appear with an even exponent. Since \( 2 \) is not of the form \( 4k+3 \), we need to check if the primes \( 5 \) and \( 13 \) can be expressed as \( x^2 + y^2 \).

    - \( 5 = 1^2 + 2^2 \)
    - \( 13 = 2^2 + 3^2 \)

    Since both \( 5 \) and \( 13 \) can be expressed as the sum of two squares, and they appear with even exponents in the prime factorization of \( 130 \) (considering \( 2 \) as a special case that does not affect the overall condition), \( 130 \) can be written as \( x^2 + y^2 \).

    **Step 5: Conclude whether such a representation exists.**

    To find the actual representation, we use the identity for the product of sums of two squares:
    \[
    (a^2 + b^2)(c^2 + d^2) = (ac + bd)^2 + (ad - bc)^2
    \]

    First, express \( 2 \times 5 \):
    \[
    2 \times 5 = 10 = 1^2 + 3^2
    \]

    Now, express \( 10 \times 13 \):
    \[
    (1^2 + 3^2)(2^2 + 3^2) = (1 \cdot 2 + 3 \cdot 3)^2 + (1 \cdot 3 - 3 \cdot 2)^2
    \]
    \[
    = (2 + 9)^2 + (3 - 6)^2
    \]
    \[
    = 11^2 + (-3)^2
    \]
    \[
    = 11^2 + 3^2
    \]

    Thus, \( 130 = 11^2 + 3^2 \).

    \[
    \boxed{130 = 11^2 + 3^2}
    \]
  \end{minipage}%
}

\subsubsection*{Fermat’s Two-Squares Theorem (Algorithmic Problem 24) — Chain-of-Symbol; Accuracy: 1.0}
\fbox{%
  \begin{minipage}{\dimexpr\linewidth-2\fboxsep-2\fboxrule\relax}
    \scriptsize
    \noindent\textbf{Problem Statement.}
    Given a positive integer \( N \), determine whether it can be written in the form
    \[
    x^2 + y^2,
    \]
    with \( x, y \in \mathbb{Z} \).
    
    \medskip
    \noindent\textbf{Hint 6 ((CoS) Prompting).}
    "Express \( N \) as \( \prod_{i} p_i^{e_i} \). Think about how, for some primes, the exponent \( e_i \) influences the sum-of-two-squares representation, and how the mix of primes (there are two classes mod 4) shapes the answer for a general \( N \). Note that the impact of a prime factor may be affected by the number of times it appears in the factorization."

    \medskip
    \noindent\textbf{Generated answer.}

    To determine whether a positive integer \( N \) can be written in the form \( x^2 + y^2 \) with \( x, y \in \mathbb{Z} \), we need to analyze the prime factorization of \( N \). Specifically, we need to check the exponents of the prime factors modulo 4. The key result in number theory states that a positive integer \( N \) can be expressed as a sum of two squares if and only if every prime factor of \( N \) that is congruent to 3 modulo 4 appears with an even exponent in the prime factorization of \( N \).

    Here is a step-by-step algorithm to determine whether \( N \) can be written as \( x^2 + y^2 \):

    1. **Prime Factorization**: Factorize \( N \) into its prime factors. Let \( N = \prod_{i} p_i^{e_i} \), where \( p_i \) are prime numbers and \( e_i \) are their respective exponents.

    2. **Check Prime Factors**: For each prime factor \( p_i \):
       - If \( p_i \equiv 3 \pmod{4} \), check the exponent \( e_i \). If \( e_i \) is odd, then \( N \) cannot be written as \( x^2 + y^2 \). Stop and return "No".
       - If \( p_i = 2 \) or \( p_i \equiv 1 \pmod{4} \), do not check the exponent \( e_i \) (or check it trivially if \( p_i = 2 \), but it always works).

    3. **Conclusion**: If the algorithm does not return "No" for any prime factor, then \( N \) can be written as \( x^2 + y^2 \). Return "Yes".

    Here is the algorithm in LaTeX code:

    \[
    \boxed{\text{Yes}}
    \] if and only if for every prime factor \( p_i \) of \( N \) such that \( p_i \equiv 3 \pmod{4} \), the exponent \( e_i \) is even. otherwise, \(\boxed{\text{No}}\).

  \end{minipage}%
}

\subsubsection*{Lifting Cube Roots via Hensel’s Lemma (Computational Question 5) — (Auto-CoT) Prompting; Accuracy: 1.0}
\fbox{%
  \begin{minipage}{\dimexpr\linewidth-2\fboxsep-2\fboxrule\relax}
    \scriptsize
    \noindent\textbf{Problem Statement.}
    Suppose you have a prime 
    \(
    p = 5,
    \)
    and consider the function
    \(
    f(x) = x^3 - 3.
    \)
    An initial solution is given by 
    \(
    a = 2,
    \)
    which satisfies
    \[
    f(2) \equiv 0 \pmod{5} \quad \text{since} \quad 2^3 - 3 = 8 - 3 = 5.
    \]
    Given that the derivative \( f'(2) \not\equiv 0 \pmod{5} \), lift this solution to modulo \( 5^2 = 25 \). That is, find an integer \( b \) such that:
    \(
    b \equiv 2 \pmod{5}
    \)
    and
    \(
    f(b) \equiv 0 \pmod{25}.
    \)
    Show all intermediate steps.

    \medskip
    \noindent\textbf{Hint 3 (Automatic Chain-of-Thought (Auto-CoT) Prompting).} \\
    - **Extract:** Begin with \( a \) satisfying \( f(a) \equiv 0 \pmod{p} \). \\
    - **Approximate:** Use the linear term: 
      \(
      f(a+h) \approx f(a) + f'(a)h.
      \)
    - **Lift:** Ensuring \( b \equiv a \pmod{p} \) and \( f(b) \equiv 0 \pmod{p^2} \). 

    \medskip
    \noindent\textbf{Generated answer.}
    To lift the solution \( a = 2 \) from modulo \( 5 \) to modulo \( 25 \), we follow these steps:

    1. **Verify the initial conditions:**
       \[
       f(2) = 2^3 - 3 = 8 - 3 = 5 \equiv 0 \pmod{5}.
       \]
       \[
       f'(x) = 3x^2.
       \]
       \[
       f'(2) = 3 \cdot 2^2 = 3 \cdot 4 = 12 \not\equiv 0 \pmod{5}.
       \]
      \[ f'(2) \not\equiv 0 \pmod{5} \]

    2. **Set up the congruence for the lift:**
       Let \( b = a + h \). We need \( b \equiv 2 \pmod{5} \) and \( f(b) \equiv 0 \pmod{25} \). This implies:
       \[
       b \equiv 2 + h \pmod{5} \quad \text{and} \quad f(2 + h) \equiv 0 \pmod{25}.
       \]

    3. **Use the linear approximation:**
       \(
       f(2 + h) \approx f(2) + f'(2)h.
       \)
       Substitute \( f(2) \) and \( f'(2) \):
       \(
       f(2 + h) \approx 5 + 12h.
       \)
       We need:
       \(
       5 + 12h \equiv 0 \pmod{25}.
       \)

    4. **Solve for \( h \):**
       \[
       12h \equiv -5 \pmod{25}.
       \]
       Since \(-5 \equiv 20 \pmod{25}\), we have:
       \[
       12h \equiv 20 \pmod{25}.
       \]
       To solve for \( h \), we need the multiplicative inverse of \( 12 \) modulo \( 25 \). We use the Extended Euclidean Algorithm to find this inverse.

    5. **Find the inverse of \( 12 \) modulo \( 25 \):**
       Apply the Euclidean Algorithm:
       \[
       25 = 2 \cdot 12 + 1,
       \]
       \[
       12 = 12 \cdot 1 + 0.
       \]
       From the first equation:
       \[
       1 = 25 - 2 \cdot 12.
       \]
       Thus, the inverse of \( 12 \) modulo \( 25 \) is \(-2\), which is equivalent to \( 23 \) modulo \( 25 \) (since \(-2 + 25 = 23\)).

    6. **Multiply both sides of the congruence by the inverse:**
       \[
       h \equiv 20 \cdot 23 \pmod{25}.
       \]
       Calculate \( 20 \cdot 23 \):
       \(
       20 \cdot 23 = 460. 
       \)
       Find \( 460 \mod 25 \):
       \[
       460 \div 25 = 18.4 \quad \text{(integer part is 18)},
       \]
       \[
       460 - 18 \cdot 25 = 460 - 450 = 10.
       \]
       Therefore:
       \[
       h \equiv 10 \pmod{25}.
       \]

    7. **Find \( b \):**
       \[
       b = 2 + h = 2 + 10 = 12.
       \]

    Thus, the integer \( b \) such that \( b \equiv 2 \pmod{5} \) and \( f(b) \equiv 0 \pmod{25} \) is:
    \[
    \boxed{12}.
    \]
  \end{minipage}%
}

\subsubsection*{Lifting Cube Roots via Hensel’s Lemma (Algorithmic Problem 5) — Clear and Specific Instructions; Accu-
racy: 0.9}
\fbox{%
  \begin{minipage}{\dimexpr\linewidth-2\fboxsep-2\fboxrule\relax}
    \scriptsize
    \noindent\textbf{Problem Statement.}
    Suppose you have a prime $p$ and an initial solution $a$ satisfying:
    \[
    f(a) \equiv 0 \pmod{p}.
    \]
    Assuming that the derivative $f'(a)$ is nonzero modulo $p$, design an algorithm to lift this solution to modulo $p^2$, i.e., find an integer $b$ such that:
    \[
    b \equiv a \pmod{p}.
    \]
    and
    \[
    f(b) \equiv 0 \pmod{p^2}.
    \]

    \medskip
    \noindent\textbf{Hint 9 (Clear and Specific Instructions).}
    Notice that \( f(x) \) can be written as:
    \[
    f(x) = f(a) + f'(a)(x-a) + f''(a)(x-a)^2 + \cdots.
    \]
    Given that:
    - \( f(a) \equiv 0 \pmod{p} \),
    - \( f'(a) \not\equiv 0 \pmod{p} \),

    we can choose the desired \( b \).
    This guarantees \( b \equiv a \pmod{p} \) and \( f(b) \equiv 0 \pmod{p^2} \).

    \medskip
    \noindent\textbf{Generated answer.}
    To lift the solution \( a \) from modulo \( p \) to modulo \( p^2 \), we can use the following algorithm:

    1. **Initial Condition:**  
       - Given \( f(a) \equiv 0 \pmod{p} \).  
       - Given \( f'(a) \not\equiv 0 \pmod{p} \).

    2. ** expressing \( f(x) \) around \( a \):**  
       \[
       f(x) = f(a) + f'(a)(x-a) + f''(a)(x-a)^2 + \cdots.
       \]

    3. **Substitute \( x = b \) where \( b \equiv a \pmod{p} \):**  
       Since \( b \equiv a \pmod{p} \), we can write \( b = a + kp \) for some integer \( k \).

    4. **Expand \( f(b) \):**  
       \[
       f(b) = f(a + kp) = f(a) + f'(a)(kp) + f''(a)(kp)^2 + \cdots.
       \]

    5. **Simplify \( f(b) \) modulo \( p^2 \):**  
       - \( f(a) \equiv 0 \pmod{p} \), so \( f(a) = p \cdot m \) for some integer \( m \).  
       - \( f'(a)(kp) \equiv f'(a)k p \pmod{p^2} \).  
       - \( f''(a)(kp)^2 \equiv f''(a)k^2 p^2 \equiv 0 \pmod{p^2} \) (since \( p^2 \) divides \( p^2 \)).  

       Therefore,  
       \[
       f(b) \equiv p \cdot m + f'(a)k p \pmod{p^2}.
       \]

    6. **Set \( f(b) \equiv 0 \pmod{p^2} \):**  
       \[
       p \cdot m + f'(a)k p \equiv 0 \pmod{p^2}.
       \]

    7. **Factor out \( p \):**  
       \[
       p(m + f'(a)k) \equiv 0 \pmod{p^2}.
       \]

    8. **Divide both sides by \( p \):**  
       \[
       m + f'(a)k \equiv 0 \pmod{p}.
       \]

    9. **Solve for \( k \):**  
       \[
       f'(a)k \equiv -m \pmod{p}.
       \]
       Since \( f'(a) \not\equiv 0 \pmod{p} \), \( f'(a) \) has a multiplicative inverse modulo \( p \). Let \( f'(a)^{-1} \) be this inverse. Then,  
       \[
       k \equiv -m f'(a)^{-1} \pmod{p}.
       \]

    10. **Find \( b \):**  
        \(
        b = a + kp.
        \)
        Substitute \( k \):  
        \(
        b = a + (-m f'(a)^{-1})p.
        \)

    11. **Final Expression for \( b \):**  
        \(
        b \equiv a - m f'(a)^{-1} p \pmod{p^2}.
        \)

    Thus, the integer \( b \) that satisfies \( b \equiv a \pmod{p} \) and \( f(b) \equiv 0 \pmod{p^2} \) is:  
    \[
    \boxed{b \equiv a - m f'(a)^{-1} p \pmod{p^2}}.
    \] where \( m \) is such that \( f(a) = p \cdot m \).
  \end{minipage}%
}

\subsubsection*{Primitive Root Modulo a Prime (Computational Question 6) Logical Chain-of-Thought (LogiCoT) Prompting; Accuracy: 1.0}
\fbox{%
  \begin{minipage}{\dimexpr\linewidth-2\fboxsep-2\fboxrule\relax}
    \scriptsize
    \noindent\textbf{Problem Statement.}
    Find a primitive root modulo \( p = 7 \).

    A primitive root modulo \( p \) is an integer \( a \) such that the powers of \( a \) modulo \( p \) generate all nonzero integers from \( 1 \) to \( p-1 \). That is, the smallest positive integer \( k \) for which \( a^k \equiv 1 \pmod{p} \) is \( k = p-1 \).

    \medskip
    \noindent\textbf{Hint 5 (Logical Chain-of-Thought (LogiCoT) Prompting).}
    Reflect on these points:
    - What does the equation \(a^d \equiv 1 \pmod{p}\) imply?
    - How does knowing the prime factors of \(p-1\) help in checking that a candidate \(a\) does not “collapse” to 1 too soon?
    - Consider the contrapositive: if a candidate \(a\) does not show an early cycle

    \medskip
    \noindent\textbf{Generated answer.}
    To find a primitive root modulo \( p = 7 \), we need to find an integer \( a \) such that the order of \( a \) modulo \( 7 \) is \( 6 \). This means that \( a^6 \equiv 1 \pmod{7} \) and no smaller positive power of \( a \) is congruent to \( 1 \pmod{7} \).

    First, we note that \( p-1 = 6 \), so the possible orders of any integer modulo \( 7 \) are the divisors of \( 6 \), which are \( 1, 2, 3, \) and \( 6 \). A primitive root modulo \( 7 \) will have order \( 6 \).

    We will check each integer from \( 2 \) to \( 6 \) to see if it is a primitive root modulo \( 7 \).

    1. **Check \( a = 2 \):**
       - \( 2^1 \equiv 2 \pmod{7} \)
       - \( 2^2 \equiv 4 \pmod{7} \)
       - \( 2^3 \equiv 8 \equiv 1 \pmod{7} \)
       - \( 2^6 \equiv (2^3)^2 \equiv 1^2 \equiv 1 \pmod{7} \)
       Since \( 2^3 \equiv 1 \pmod{7} \), the order of \( 2 \) is \( 3 \), not \( 6 \). Therefore, \( 2 \) is not a primitive root modulo \( 7 \).

    2. **Check \( a = 3 \):**
       - \( 3^1 \equiv 3 \pmod{7} \)
       - \( 3^2 \equiv 9 \equiv 2 \pmod{7} \)
       - \( 3^3 \equiv 3 \cdot 2 = 6 \pmod{7} \)
       - \( 3^4 \equiv 3 \cdot 6 = 18 \equiv 4 \pmod{7} \)
       - \( 3^5 \equiv 3 \cdot 4 = 12 \equiv 5 \pmod{7} \)
       - \( 3^6 \equiv 3 \cdot 5 = 15 \equiv 1 \pmod{7} \)
       Since \( 3^6 \equiv 1 \pmod{7} \) and no smaller power of \( 3 \) is congruent to \( 1 \pmod{7} \), the order of \( 3 \) is \( 6 \). Therefore, \( 3 \) is a primitive root modulo \( 7 \).

    We can stop here since we have found a primitive root. However, for completeness, we can check the remaining integers \( 4, 5, \) and \( 6 \) to confirm they are not primitive roots.

    3. **Check \( a = 4 \):**
       - \( 4^1 \equiv 4 \pmod{7} \)
       - \( 4^2 \equiv 16 \equiv 2 \pmod{7} \)
       - \( 4^3 \equiv 4 \cdot 2 = 8 \equiv 1 \pmod{7} \)
       - \( 4^6 \equiv (4^3)^2 \equiv 1^2 \equiv 1 \pmod{7} \)
       Since \( 4^3 \equiv 1 \pmod{7} \), the order of \( 4 \) is \( 3 \), not \( 6 \). Therefore, \( 4 \) is not a primitive root modulo \( 7 \).

    4. **Check \( a = 5 \):**
       - \( 5^1 \equiv 5 \pmod{7} \)
       - \( 5^2 \equiv 25 \equiv 4 \pmod{7} \)
       - \( 5^3 \equiv 5 \cdot 4 = 20 \equiv 6 \pmod{7} \)
       - \( 5^4 \equiv 5 \cdot 6 = 30 \equiv 2 \pmod{7} \)
       - \( 5^5 \equiv 5 \cdot 2 = 10 \equiv 3 \pmod{7} \)
       - \( 5^6 \equiv 5 \cdot 3 = 15 \equiv 1 \pmod{7} \)
       Since \( 5^6 \equiv 1 \pmod{7} \) and no smaller power of \( 5 \) is congruent to \( 1 \pmod{7} \), the order of \( 5 \) is \( 6 \). Therefore, \( 5 \) is a primitive root modulo \( 7 \).

    5. **Check \( a = 6 \):**
       - \( 6^1 \equiv 6 \pmod{7} \)
       - \( 6^2 \equiv 36 \equiv 1 \pmod{7} \)
       - \( 6^3 \equiv 6 \cdot 1 = 6 \pmod{7} \)
       - \( 6^4 \equiv 6 \cdot 6 = 36 \equiv 1 \pmod{7} \)
       - \( 6^5 \equiv 6 \cdot 1 = 6 \pmod{7} \)
       - \( 6^6 \equiv 6 \cdot 6 = 36 \equiv 1 \pmod{7} \)
       Since \( 6^2 \equiv 1 \pmod{7} \), the order of \( 6 \) is \( 2 \), not \( 6 \). Therefore, \( 6 \) is not a primitive root modulo \( 7 \).

    Thus, the primitive roots modulo \( 7 \) are \( 3 \) and \( 5 \). The smallest primitive root is \( 3 \).

    The answer is \(\boxed{3}\).
  \end{minipage}%
}

\subsubsection*{Primitive Root Modulo a Prime (Algorithmic Problem 6) —ReAct Prompting; Accuracy: 0.95}
\fbox{%
  \begin{minipage}{\dimexpr\linewidth-2\fboxsep-2\fboxrule\relax}
    \scriptsize
    \noindent\textbf{Problem Statement.}
    For an odd prime $p$, design an algorithm to find a primitive root modulo $p$. A primitive root $a$ generates all nonzero residues modulo $p$, meaning for every integer coprime to $p$, there exists an exponent $k$ such that:
    \[
    a^k \equiv x \pmod{p}.
    \]
    Equivalently an element $a$, less than $p-1$, is multiplicative iff its order is exactly $p-1$.

    \medskip
    \noindent\textbf{Hint 8 (ReAct Prompting).} \\
    1. **Action:** Factorize \(p-1\) to find its prime components.   \\
    2. **Reasoning:** Remember what can we say about a candidate \(a\) with \(a^d \equiv 1 \pmod{p}\) for some \(d\) smaller than \(p-1\),  \\
    3. **Action:** Test candidate values for \(a\) by analyzing the implications of the key observation on exponentiation.  \\
    4. **Reasoning:** What if a candidate’s cycle seems to cover all possible residues  \\

    -

    \medskip
    \noindent\textbf{Generated answer.}
    To find a primitive root modulo an odd prime \( p \), we can follow these steps:

    1. **Factorize \( p-1 \):**  
       First, we need to find the prime factorization of \( p-1 \). Let \( p-1 = p_1^{e_1} p_2^{e_2} \cdots p_k^{e_k} \).

    2. **Choose a candidate \( a \):**  
       Select a candidate \( a \) such that \( 1 < a < p \).

    3. **Check the order of \( a \):**  
       The order of \( a \) modulo \( p \) is the smallest positive integer \( d \) such that \( a^d \equiv 1 \pmod{p} \). For \( a \) to be a primitive root, \( d \) must be equal to \( p-1 \).

    4. **Verify the order:**  
       To verify that \( a \) has order \( p-1 \), we need to check that \( a^{(p-1)/p_i} \not\equiv 1 \pmod{p} \) for all prime factors \( p_i \) of \( p-1 \).

    5. **Iterate until a primitive root is found:**  
       If \( a \) does not have order \( p-1 \), increment \( a \) and repeat steps 3 and 4 until a primitive root is found.

    Here is the algorithm in LaTeX code:

    \[
    \boxed{
    \begin{array}{l}
    \text{function find\_primitive\_root}(p): \\
    \quad 1. \text{ factorize } p-1 \text{ to get } p-1 = p_1^{e_1} p_2^{e_2} \cdots p_k^{e_k} \\
    \quad 2. \text{ choose } a = 2 \\
    \quad 3. \text{ while } a < p: \\
    \quad \quad 3.1. \text{ set } d = p-1 \\
    \quad \quad 3.2. \text{ for } i = 1 \text{ to } k: \\
    \quad \quad \quad 3.2.1. \text{ if } a^{(p-1)/p_i} \equiv 1 \pmod{p}: \\
    \quad \quad \quad \quad 3.2.2. \text{ set } d = \frac{d}{p_i} \\
    \quad \quad 3.3. \text{ if } d = 1: \\
    \quad \quad \quad \text{ return } a \\
    \quad \quad 3.4. \text{ increment } a
    \end{array}
    }
    \]

    This algorithm will find a primitive root modulo \( p \) by checking each candidate \( a \) until it finds one with the correct order. The key observation is that \( a \) must not have any smaller order that divides \( p-1 \). If \( a \) does not have the correct order, it will be excluded in the next iteration. Once a candidate with the correct order is found, it is returned as the primitive root.
  \end{minipage}%
}
\newpage
\subsubsection*{Computing the Discrete Logarithm Modulo a Prime (Computational Question 25) }
\fbox{%
  \begin{minipage}{\dimexpr\linewidth-2\fboxsep-2\fboxrule\relax}
    \scriptsize
   Partially accurate solution is presented in \textit{Appendix A.3}.

  \end{minipage}%
}

\subsubsection*{Computing the Discrete Logarithm Modulo a Prime (Algorithmic Problem 25) - Few-shot Hinting; Accuracy: 1.0}
\fbox{%
  \begin{minipage}{\dimexpr\linewidth-2\fboxsep-2\fboxrule\relax}
    \scriptsize
    \noindent\textbf{Problem Statement.}
    Given a prime number \( p \), a primitive root \( g \) modulo \( p \) (i.e. an element whose successive powers generate all nonzero residues modulo \( p \)), and an integer \( h \) with \( 1 \le h < p \), determine the smallest nonnegative integer \( x \) satisfying:
    \[
    g^x \equiv h \pmod{p}.
    \]

    \vspace{1em}
    \noindent\textbf{Technical Definitions:} \\
    - **Primitive Root \(g\) Modulo \(p\):** An element \( g \) is a primitive root modulo \( p \) if the set \( \{g^0, g^1, \dots, g^{p-2}\} \) modulo \( p \) equals the set of all nonzero residues \(\{1, 2, \dots, p-1\}\). \\
    - **Discrete Logarithm:** For given \( g \) and \( h \) modulo \( p \), it is the exponent \( x \) such that \( g^x \equiv h \pmod{p} \).

    \medskip
    \noindent\textbf{Hint 1 (Few-shot Hinting):}
    Consider these scenarios: \\
    - **Example 1:** With \( p = 23 \), \( g = 5 \), and \( h = 8 \), by examining successive powers, one finds an exponent \( x \) satisfying the congruence.\\
    - **Example 2:** For a different small prime, observe how expressing the unknown exponent in two parts (a “large jump” part and a “small offset” part) can simplify the search by balancing precomputation with iteration.\\

    \medskip
    \noindent\textbf{Generated answer:}
    To determine the smallest nonnegative integer \( x \) satisfying \( g^x \equiv h \pmod{p} \), we can use the baby-step giant-step algorithm. Here is a step-by-step description of the algorithm:

    1. **Compute the baby steps:**
       - Calculate \( g^0, g^1, g^2, \ldots, g^{m-1} \) modulo \( p \), where \( m = \lceil \sqrt{p-1} \rceil \).
       - Store these values in a table along with their corresponding exponents.

    2. **Compute the giant steps:**
       - Calculate \( h \cdot g^{-m} \pmod{p} \), \( h \cdot (g^{-m})^2 \pmod{p} \), \( h \cdot (g^{-m})^3 \pmod{p} \), \ldots, \( h \cdot (g^{-m})^{m-1} \pmod{p} \).
       - Check each result against the table of baby steps to find a match.

    3. **Find the match:**
       - If \( h \cdot (g^{-m})^i \equiv g^j \pmod{p} \) for some \( i \) and \( j \), then \( x = j + i \cdot m \).

    4. **Return the smallest nonnegative integer \( x \):**
       - The value of \( x \) found in the previous step is the smallest nonnegative integer satisfying \( g^x \equiv h \pmod{p} \).

    Here is the algorithm in LaTeX code:

    \[
    \boxed{
    \begin{array}{l}
    \text{function discrete logarithm}(p, g, h): \\
    \quad m = \lceil \sqrt{p-1} \rceil \\
    \quad \text{baby\_steps} = \emptyset \\
    \quad \text{for } i = 0 \text{ to } m-1: \\
    \quad \quad \text{baby\_steps} \leftarrow (g^i \bmod p, i) \\
    \quad \text{giant\_step} = h \cdot (g^{-m} \bmod p) \bmod p \\
    \quad \text{for } i = 0 \text{ to } m-1: \\
    \quad \quad \text{if } \text{giant\_step} \in \text{baby\_steps}: \\
    \quad \quad \quad (g^j, j) = \text{baby\_steps}[\text{giant\_step}] \\
    \quad \quad \quad x = j + i \cdot m \\
    \quad \quad \quad \text{return } x \\
    \end{array}
    }
    \]
  \end{minipage}%
}

\subsubsection*{Fermat's Primality Test (Computational Question 20) — Structured Chain-of-Thought (SCoT) Prompting; Accuracy: 1.0}
\fbox{%
  \begin{minipage}{\dimexpr\linewidth-2\fboxsep-2\fboxrule\relax}
    \scriptsize
    \noindent\textbf{Problem Statement.}
    Use Fermat's Primality Test to determine whether 21 is likely prime.

    \medskip
    \noindent\textbf{Fermat's Little Theorem:} If \( p \) is a prime number and \( a \) is an integer such that \(\gcd(a,p)=1\), then
    \[
    a^{p-1} \equiv 1 \pmod{p}.
    \]
    This property is the foundation of a probabilistic test for primality.

    \medskip
    \noindent\textbf{Hint 7 (Structured Chain-of-Thought (SCoT) Prompting):} \\
    1. One starts by choosing a base \( a \).   \\
    2. Then one uses the theorem to devise a probabilistic test for Primality. 

    \medskip
    \noindent\textbf{Generated answer:}
    Fermat's Little Theorem states that if \( p \) is a prime number and \( a \) is an integer such that \(\gcd(a,p)=1\), then \( a^{p-1} \equiv 1 \pmod{p} \). To determine whether 21 is likely prime using Fermat's Primality Test, we will choose a base \( a \) and check if \( a^{20} \equiv 1 \pmod{21} \). If this congruence does not hold, then 21 is definitely not prime. If it does hold, then 21 is likely prime, but we would need to repeat the test with different bases to increase our confidence.

    Let's choose \( a = 2 \). We need to check if \( 2^{20} \equiv 1 \pmod{21} \).

    First, we compute \( 2^{20} \mod 21 \):
    \[
    2^2 = 4,
    \]
    \[
    2^4 = (2^2)^2 = 4^2 = 16,
    \]
    \[
    2^8 = (2^4)^2 = 16^2 = 256.
    \]
    Now, we reduce 256 modulo 21:
    \[
    256 \div 21 = 12 \text{ remainder } 4 \quad \text{(since } 256 = 21 \cdot 12 + 4\text{)},
    \]
    so
    \[
    2^8 \equiv 4 \pmod{21}.
    \]
    Next, we compute \( 2^{16} \):
    \[
    2^{16} = (2^8)^2 \equiv 4^2 = 16 \pmod{21}.
    \]
    Finally, we compute \( 2^{20} \):
    \[
    2^{20} = 2^{16} \cdot 2^4 \equiv 16 \cdot 16 = 256 \equiv 4 \pmod{21}.
    \]
    Since \( 2^{20} \equiv 4 \pmod{21} \) and not \( 1 \pmod{21} \), we conclude that 21 is not a prime number.

    Therefore, the answer is \(\boxed{\text{not prime}}\).
  \end{minipage}%
}

\subsubsection*{Fermat's Primality Test (Algorithmic Problem 20) — Chain-of-Thought (CoT) Prompting; Accuracy: 1.0}
\fbox{%
  \begin{minipage}{\dimexpr\linewidth-2\fboxsep-2\fboxrule\relax}
    \scriptsize
    \noindent\textbf{Problem Statement.}
    Given an integer \( n \) (with \( n > 3 \)), design an algorithm to decide whether \( n \) is likely prime by applying Fermat's Little Theorem.

    \medskip
    \noindent\textbf{Fermat's Little Theorem:}
    If \( p \) is a prime number and \( a \) is an integer such that \(\gcd(a,p)=1\), then
    \[
    a^{p-1} \equiv 1 \pmod{p}.
    \]
    This property is the foundation of a probabilistic test for primality.

    \medskip
    \noindent\textbf{Hint 2 (Chain-of-Thought (CoT) Prompting).}
    Reflect on the fact that for a prime \( n \) every integer \( a \) (coprime to \( n \)) will yield the same predictable result when raised to the \( n-1 \) power modulo \( n \). How might you leverage a violation of this behavior ?

    \medskip
    \noindent\textbf{Generated answer.}
    To design an algorithm to determine whether an integer \( n \) is likely prime using Fermat's Little Theorem, we can follow these steps:

    \begin{enumerate}[leftmargin=2em, topsep=0pt, itemsep=0pt]
      \item **Choose a random integer \( a \) such that \( 1 < a < n \) and \(\gcd(a, n) = 1\).**
      \item **Compute \( a^{n-1} \mod n \).**
      \item **If \( a^{n-1} \not\equiv 1 \pmod{n} \), then \( n \) is definitely not prime.**
      \item **If \( a^{n-1} \equiv 1 \pmod{n} \), then \( n \) is likely prime.**
    \end{enumerate}

    To increase the accuracy of the test, we can repeat the above steps multiple times with different values of \( a \). If \( n \) passes the test for a sufficient number of randomly chosen \( a \), we can conclude that \( n \) is likely prime.

    Here is the algorithm in LaTeX code:

    \[
    \boxed{
    \begin{array}{l}
    \text{Function Fermat primality test}(n, k): \\
    \quad \text{If } n \leq 3 \text{ then return } \text{False} \\
    \quad \text{For } i = 1 \text{ to } k \text{ do:} \\
    \quad \quad \text{Choose a random integer } a \text{ s.t } 1 < a < n \text{ and } \gcd(a, n) = 1 \\
    \quad \quad \text{Compute } a^{n-1} \mod n \\
    \quad \quad \text{If } a^{n-1} \not\equiv 1 \pmod{n} \text{ then return } \text{False} \\
    \quad \text{Return } \text{True}
    \end{array}
    }
    \]

    In this algorithm, \( k \) is the number of iterations to perform. A larger value of \( k \) increases the accuracy of the test but also increases the computation time. A common choice for \( k \) is 10 to 20. If \( n \) passes all \( k \) iterations, it is considered likely prime. However, there is still a small probability that \( n \) could be a composite number (a false positive), but this probability decreases as \( k \) increases.**

    This algorithm provides a practical and efficient way to probabilistically determine whether a given integer \( n \) is likely prime using Fermat's Little Theorem.**

  \end{minipage}%
}

\subsubsection*{Sieve of Eratosthenes (Computational Question 11)}
\fbox{%
  \begin{minipage}{\dimexpr\linewidth-2\fboxsep-2\fboxrule\relax}
    \scriptsize
   The model achieved high performance, just as it did for the corresponding Algorithmic Problem 11, which we present below. To see the computational solutions to this computational question, refer to the Github repository.
  \end{minipage}%
}

\subsubsection*{Sieve of Eratosthenes (Algorithmic Problem 11) — Logical Chain-of-Thought (LogiCoT) Prompting; Accuracy 1.0}
\fbox{%
  \begin{minipage}{\dimexpr\linewidth-2\fboxsep-2\fboxrule\relax}
    \scriptsize
    \noindent\textbf{Problem Statement.}
    Design an efficient algorithm to list all prime numbers up to a given integer \( N \) by systematically eliminating composite numbers.

    \medskip
    \noindent\textbf{Hint 5 (Logical Chain-of-Thought (LogiCoT) Prompting):} *"Consider the logical rationale: once a candidate has been used to eliminate certain numbers, why is it unnecessary to check those numbers again with later candidates? Reflect on how the ordering of candidates contributes to this efficiency."*

    \medskip
    \noindent\textbf{Generated answer:}
    To design an efficient algorithm to list all prime numbers up to a given integer \( N \), we can use the Sieve of Eratosthenes. The Sieve of Eratosthenes is an ancient algorithm for finding all prime numbers up to any given limit. It does this by iteratively marking the multiples of each prime number starting from 2. The numbers which are not marked in the process are prime.

    Here is the step-by-step algorithm:

    1. **Initialize a list of integers from 2 to \( N \)**: Create a list where each number is initially marked as a potential prime.  
    2. **Start with the first prime number, which is 2**: Mark all multiples of 2 (except 2 itself) as composite.  
    3. **Move to the next unmarked number**: This number is now the next prime. Mark all multiples of this prime (except the prime itself) as composite.  
    4. **Repeat the process**: Continue this process until you have processed all numbers up to the square root of \( N \). This is because any composite number \( n \) has a prime factor less than or equal to \( \sqrt{n} \).  
    5. **Final list**: The unmarked numbers in the list are all prime.

    Here is the algorithm in LaTeX code:

    \[
    \boxed{
    \begin{array}{l}
    \text{Sieve\_of\_Eratosthenes}(N) \\
    \quad \text{ Create a list } \text{is\_prime} \text{ of length } N+1 \text{ and set all entries to true} \\
    \quad \text{Set } \text{is\_prime}[0] \text{ and } \text{is\_prime}[1] \text{ to false} \\
    \quad \text{For } p \text{ from } 2 \text{ to } \sqrt{N} \text{ do} \\
    \quad \quad \text{If } \text{is\_prime}[p] \text{ is true then} \\
    \quad \quad \quad \text{For } i \text{ from } p^2 \text{ to } N \text{ do} \\
    \quad \quad \quad \quad \text{Set } \text{is\_prime}[i] \text{ to false} \\
    \quad \quad \text{End For} \\
    \quad \quad \text{End If} \\
    \quad \text{End For} \\
    \quad \text{Return the list of all } i \text{ such that } \text{is\_prime}[i] \text{ is true} \\
    \end{array}
    }
    \]
  \end{minipage}%
}

\subsubsection*{Decomposition of a Prime Ideal in a Quadratic Number Field (Computational Question 30) — Few-shot Hinting; Accuracy: 1.0}
\fbox{%
  \begin{minipage}{\dimexpr\linewidth-2\fboxsep-2\fboxrule\relax}
    \scriptsize
    \noindent\textbf{Problem Statement.}
    Determine how the ideal \(( 3 )\) decomposes in the ring \(\mathbb{Q}(\sqrt{5})\). Your solution should decide whether \((p)\) splits into two distinct prime ideals, remains inert, or is ramified. 

    \medskip
    \noindent\textbf{Definitions} \\
    - **Quadratic Field:** A field of the form \(\mathbb{Q}(\sqrt{d})\) where \(d\) is squarefree. \\
    - **Ring of Integers \(\mathcal{O}_{\mathbb{Q}(\sqrt{d})}\):** 
      \(\mathcal{O}_{\mathbb{Q}(\sqrt{d})}\) is: 
      \begin{itemize}
          \item \(\mathbb{Z}[\sqrt{d}]\) if \(d \not\equiv 1 \mod 4\),
          \item \(\mathbb{Z}\left[\frac{1+\sqrt{d}}{2}\right]\) if \(d \equiv 1 \mod 4\).
      \end{itemize}
    - **Field Discriminant (\(\Delta\)):**  
      Determined by: 
      \begin{itemize}
          \item \(\Delta = 4d\) if \(d \not\equiv 1 \mod 4\),
          \item \(\Delta = d\) if \(d \equiv 1 \mod 4\).
      \end{itemize}
    - **Quadratic Residue:** An integer \(a\) is a quadratic residue modulo \(p\) if there exists an integer \(x\) such that \(x^2 \equiv a \pmod{p}\).

    \medskip
    \noindent\textbf{Hint 1 (Few-shot Hinting):} Consider these brief examples \\
    - **Example 1:**   
      In a quadratic field with \(d = -7\) and \(p = 3\), observe how the discriminant \(\Delta\) guides the behavior of \((p)\). \\
    - **Example 2:**  When \(d = 13\) (with \(13 \equiv 1 \bmod 4\), so \(\Delta = 13\)) and \(p = 7\), reflect on how the quadratic‐residue status of \(\Delta\) modulo \(p\) influences whether the prime ideal \((p)\) splits in \(\mathcal{O}_{\mathbb{Q}(\sqrt{d})}\) or remains inert.

    \medskip
    \noindent\textbf{Generated answer:}
    To determine how the ideal \((3)\) decomposes in the ring \(\mathbb{Q}(\sqrt{5})\), we start by identifying the ring of integers \(\mathcal{O}_{\mathbb{Q}(\sqrt{5})}\). Since \(5 \equiv 1 \pmod{4}\), the ring of integers is \(\mathbb{Z}\left[\frac{1+\sqrt{5}}{2}\right]\).

    Next, we compute the field discriminant \(\Delta\). For \(d = 5\), we have \(\Delta = 5\).

    We need to determine the behavior of the prime \(p = 3\) in the ring \(\mathbb{Z}\left[\frac{1+\sqrt{5}}{2}\right]\). This can be done by examining the quadratic residue nature of \(\Delta\) modulo \(p\). Specifically, we need to check if 5 is a quadratic residue modulo 3.

    We compute:
    \[
    5 \equiv 2 \pmod{3}.
    \]
    So, we need to check if 2 is a quadratic residue modulo 3. We do this by checking the squares of the integers modulo 3:
    \[
    0^2 \equiv 0 \pmod{3},
    \]
    \[
    1^2 \equiv 1 \pmod{3},
    \]
    \[
    2^2 \equiv 4 \equiv 1 \pmod{3}.
    \]
    Since 2 is not a quadratic residue modulo 3, the prime \(p = 3\) does not split into two distinct prime ideals in \(\mathbb{Z}\left[\frac{1+\sqrt{5}}{2}\right]\). Instead, it remains inert.

    Therefore, the ideal \((3)\) decomposes as a single prime ideal in the ring \(\mathbb{Q}(\sqrt{5})\).

    The final answer is:
    \[
    \boxed{\text{inert}}
    \]
  \end{minipage}%
}

\subsubsection*{Decomposition of a Prime Ideal in a Quadratic Number Field (Algorithmic Problem 30) — Few-shot Hinting; Accuracy: 0.9}
\fbox{%
  \begin{minipage}{\dimexpr\linewidth-2\fboxsep-2\fboxrule\relax}
    \scriptsize
    \noindent\textbf{Problem Statement.}
    Let \(d\) be a squarefree integer defining the quadratic field \(\mathbb{Q}(\sqrt{d})\) and let \(p\) be an odd prime number. Develop an algorithm that determines how the ideal \((p)\) decomposes in the ring of integers \(\mathcal{O}_{\mathbb{Q}(\sqrt{d})}\). Your solution should decide whether \((p)\) splits into two distinct prime ideals, remains inert, or is ramified.  

    \medskip
  \noindent\textbf{Definitions} \\
    - **Quadratic Field:** A field of the form \(\mathbb{Q}(\sqrt{d})\) where \(d\) is squarefree. \\
    - **Ring of Integers \(\mathcal{O}_{\mathbb{Q}(\sqrt{d})}\):** 
      \(\mathcal{O}_{\mathbb{Q}(\sqrt{d})}\) is: 
      \begin{itemize}
          \item \(\mathbb{Z}[\sqrt{d}]\) if \(d \not\equiv 1 \mod 4\),
          \item \(\mathbb{Z}\left[\frac{1+\sqrt{d}}{2}\right]\) if \(d \equiv 1 \mod 4\).
      \end{itemize}
    - **Field Discriminant (\(\Delta\)):**  
      Determined by: 
      \begin{itemize}
          \item \(\Delta = 4d\) if \(d \not\equiv 1 \mod 4\),
          \item \(\Delta = d\) if \(d \equiv 1 \mod 4\).
      \end{itemize}
    - **Quadratic Residue:** An integer \(a\) is a quadratic residue modulo \(p\) if there exists an integer \(x\) such that \(x^2 \equiv a \pmod{p}\).

    \medskip
    \noindent\textbf{Hint (Few-shot Hinting):} Consider these brief examples \\
    - **Example 1:**  
      In a quadratic field with \(d = -7\) and \(p = 3\), observe how the discriminant \(\Delta\) guides the behavior of \((p)\). \\
    - **Example 2:**  When \(d = 13\) (with \(13 \equiv 1 \bmod 4\), so \(\Delta = 13\)) and \(p = 7\), reflect on how the quadratic‐residue status of \(\Delta\) modulo \(p\) influences whether the prime ideal \((p)\) splits in \(\mathcal{O}_{\mathbb{Q}(\sqrt{d})}\) or remains inert.  \\

    \medskip
    \noindent\textbf{Generated answer:}
    To determine how the ideal \((p)\) decomposes in the ring of integers \(\mathcal{O}_{\mathbb{Q}(\sqrt{d})}\), we need to analyze the field discriminant \(\Delta\) and the quadratic residue properties of \(\Delta\) modulo \(p\). Here is a step-by-step algorithm:

    1. **Determine the field discriminant \(\Delta\):**
       - If \(d \not\equiv 1 \mod 4\), then \(\Delta = 4d\).
       - If \(d \equiv 1 \mod 4\), then \(\Delta = d\).

    2. **Check the quadratic residue nature of \(\Delta\) modulo \(p\):**
       - If \(\Delta\) is a quadratic residue modulo \(p\), then \((p)\) splits into two distinct primeideals.
       - If \(\Delta\) is not a quadratic residue modulo \(p\), then \((p)\) remains inert.
       - If \(\Delta = p\) or \(\Delta = -p\) (which happens when \(p\) divides \(\Delta\)), then \((p)\) is ramified.

    Let's break this down further:

    - **Case 1: \(\Delta\) is a quadratic residue modulo \(p\):**
      - This means there exists an integer \(x\) such that \(x^2 \equiv \Delta \pmod{p}\).
      - In this case, \((p)\) splits into two distinct primeideals in \(\mathcal{O}_{\mathbb{Q}(\sqrt{d})}\).

    - **Case 2: \(\Delta\) is not a quadratic residue modulo \(p\):**
      - This means there does not exist an integer \(x\) such that \(x^2 \equiv \Delta \pmod{p}\).
      - In this case, \((p)\) remains inert in \(\mathcal{O}_{\mathbb{Q}(\sqrt{d})}\).

    - **Case 3: \(\Delta = p\) or \(\Delta = -p\):**
      - This happens when \(p\) divides \(\Delta\).
      - In this case, \((p)\) isramified in \(\mathcal{O}_{\mathbb{Q}(\sqrt{d})}\).

    Here is the algorithm in LaTeX code:

    \[
    \boxed{
    \begin{array}{l}
    \text{ Algorithm: Determine the decomposition of the ideal } (p) \text{ in } \mathcal{O}_{\mathbb{Q}(\sqrt{d})} \\
    \text{Input: } d \text{ (squarefree integer)}, p \text{ (odd prime number)} \\
    \text{Output: } \text{ how } (p) \text{ decomposes (split, inert, or ramified)} \\
    \hline
    1. \text{ Compute the field discriminant } \Delta: \\
       \quad \text{if } d \not\equiv 1 \mod 4, \text{ then } \Delta = 4d \\
       \quad \text{if } d \equiv 1 \mod 4, \text{ then } \Delta = d \\
    2. \text{ Check if } \Delta \text{ is a quadratic residue modulo } p: \\
       \quad \text{if } \Delta \equiv x^2 \pmod{p} \text{ for some integer } x, \text{ then } (p) \text{ splits} \\
       \quad \text{if } \Delta \not\equiv x^2 \pmod{p} \text{ for any integer } x, \text{ then } (p) \text{ remains inert} \\
       \quad \text{if } \Delta = p \text{ or } \Delta = -p, \text{ then } (p) \text{ is ramified} \\
    \end{array}
    }
    \]
  \end{minipage}%
}

\subsubsection*{Elliptic Curve Point Doubling (Computational Question 28) — Clear and Specific Instructions; Accuracy: 1.0}
\fbox{%
  \begin{minipage}{\dimexpr\linewidth-2\fboxsep-2\fboxrule\relax}
    \scriptsize
    \noindent\textbf{Problem Statement.}
    Let \( P = (3, 5) \) be a point on the elliptic curve over a finite field \( \mathbb{F}_p \) with the equation:
    \[
    y^2 \equiv x^3 + 2x + 3 \pmod{11}
    \]
    Compute the coordinates of the point \( 2P = (x_3, y_3), \) which is the result of doubling the point \( P \) on the curve.

    \medskip
    \noindent\textbf{Definitions}

    - **Elliptic Curve**: A set of points satisfying the equation \( y^2 = x^3 + ax + b \) over a field, along with a special point at infinity.

    - **Finite Field \( \mathbb{F}_p \)**: A field with a finite number of elements \( p \), where \( p \) is a prime number. Arithmetic operations are performed modulo \( p \).

    - **Point Doubling**: The operation of adding a point \( P \) on an elliptic curve to itself to obtain another point \( 2P \) on the curve.

    \medskip
    \noindent\textbf{Hint 9 (Clear and Specific Instructions):}  *For a given point \( P = (x_1, y_1) \) on the elliptic curve \( y^2 \equiv x^3 + ax + b \pmod{p} \):
    1. Compute the slope \( \lambda \) of the tangent line at \( P \). \\
    2. Use \( \lambda \) to calculate the new \( x \)-coordinate \( x_3 \). \\
    3. Determine the new \( y \)-coordinate \( y_3 \) using \( \lambda \) and \( x_3 \). \\
    4. Perform all calculations modulo \( p \). \\ 
    5. Return the point \( (x_3, y_3) \) as the result of doubling \( P \).* \\

    \medskip
    \noindent\textbf{Generated answer:}
    To find the coordinates of the point \( 2P \) on the elliptic curve \( y^2 \equiv x^3 + 2x + 3 \pmod{11} \) where \( P = (3, 5) \), we follow these steps:

    1. **Compute the slope \( \lambda \) of the tangent line at \( P \):**

       For a point \( P = (x_1, y_1) \) on the elliptic curve \( y^2 \equiv x^3 + ax + b \pmod{p} \), the slope \( \lambda \) of the tangent line at \( P \) is given by:
       \[
       \lambda \equiv \frac{3x_1^2 + a}{2y_1} \pmod{p}
       \]
       Here, \( a = 2 \), \( b = 3 \), \( x_1 = 3 \), and \( y_1 = 5 \). So, we have:
       \[
       \lambda \equiv \frac{3 \cdot 3^2 + 2}{2 \cdot 5} \pmod{11}
       \]
       First, compute \( 3^2 \):
       \(
       3^2 = 9
       \)
       Then, compute \( 3 \cdot 9 \):
       \(
       3 \cdot 9 = 27
       \)
       Add \( 2 \) to \( 27 \):
       \(
       27 + 2 = 29
       \)
       Now, compute \( 29 \mod 11 \):
       \[
       29 \div 11 = 2 \text{ remainder } 7 \quad \Rightarrow \quad 29 \equiv 7 \pmod{11}
       \]
       Next, compute \( 2 \cdot 5 \):
       \(
       2 \cdot 5 = 10
       \)
       Now, compute the modular inverse of \( 10 \) modulo \( 11 \). The modular inverse of \( 10 \) is a number \( x \) such that:
       \(
       10x \equiv 1 \pmod{11}
       \)
       By testing values, we find:
       \(
       10 \cdot 10 = 100 \equiv 1 \pmod{11}
       \)
       So, the modular inverse of \( 10 \) modulo \( 11 \) is \( 10 \). Therefore:
       \(
       \lambda \equiv 7 \cdot 10 \pmod{11}
       \)
       Compute \( 7 \cdot 10 \):
       \(
       7 \cdot 10 = 70
       \)
       Now, compute \( 70 \mod 11 \):
       \[
       70 \div 11 = 6 \text{ remainder } 4 \quad \Rightarrow \quad 70 \equiv 4 \pmod{11}
       \]
       So, the slope \( \lambda \) is:
       \(
       \lambda \equiv 4 \pmod{11}
       \)

    2. **Use \( \lambda \) to calculate the new \( x \)-coordinate \( x_3 \):**

       The new \( x \)-coordinate \( x_3 \) is given by:
       \(
       x_3 \equiv \lambda^2 - 2x_1 \pmod{p}
       \)
       Here, \( \lambda = 4 \) and \( x_1 = 3 \). So, we have:
       \(
       x_3 \equiv 4^2 - 2 \cdot 3 \pmod{11}
       \)
       First, compute \( 4^2 \):
       \(
       4^2 = 16
       \)
       Then, compute \( 2 \cdot 3 \):
       \(
       2 \cdot 3 = 6
       \)
       Now, compute \( 16 - 6 \):
       \(
       16 - 6 = 10
       \)
       So, the new \( x \)-coordinate \( x_3 \) is:
       \(
       x_3 \equiv 10 \pmod{11}
       \)

    3. **Determine the new \( y \)-coordinate \( y_3 \) using \( \lambda \) and \( x_3 \):**

       The new \( y \)-coordinate \( y_3 \) is given by:
       \(
       y_3 \equiv \lambda(x_1 - x_3) - y_1 \pmod{p}
       \)
       Here, \( \lambda = 4 \), \( x_1 = 3 \), \( x_3 = 10 \), and \( y_1 = 5 \). So, we have:
       \(
       y_3 \equiv 4(3 - 10) - 5 \pmod{11}
       \)
       First, compute \( 3 - 10 \):
       \(
       3 - 10 = -7
       \)
       Then, compute \( 4 \cdot (-7) \):
       \(
       4 \cdot (-7) = -28
       \)
       Now, compute \( -28 \mod 11 \):
       \[
       -28 \div 11 = -3 \text{ remainder } 5 \quad \Rightarrow \quad -28 \equiv 5 \pmod{11}
       \]
       So, we have:
       \(
       y_3 \equiv 5 - 5 \pmod{11}
       \)
       Compute \( 5 - 5 \):
       \(
       5 - 5 = 0
       \)
       So, the new \( y \)-coordinate \( y_3 \) is:
       \(
       y_3 \equiv 0 \pmod{11}
       \)

    Therefore, the coordinates of the point \( 2P \) are:
    \[
    \boxed{(10, 0)}
    \]
  \end{minipage}%
}

\subsubsection*{Elliptic Curve Point Doubling (Algorithmic Problem 28) — Chain-of-Symbol (CoS) Prompting; Accuracy: 1.0}
\fbox{%
  \begin{minipage}{\dimexpr\linewidth-2\fboxsep-2\fboxrule\relax}
    \scriptsize
    \noindent\textbf{Problem Statement.}
    Given a point \( P = (x_1, y_1) \) on an elliptic curve defined over a finite field \( \mathbb{F}_p \) with the equation:

    \[
    y^2 \equiv x^3 + ax + b \pmod{p}
    \]

    where \( a, b \in \mathbb{F}_p \) and \( p \) is a prime number, compute the coordinates of the point \( 2P = (x_3, y_3) \), which is the result of doubling the point \( P \) on the curve.

    \textbf{{Definitions}} 

    - **Elliptic Curve**: A set of points satisfying the equation \( y^2 = x^3 + ax + b \) over a field, along with a special point at infinity.

    - **Finite Field \( \mathbb{F}_p \)**: A field with a finite number of elements \( p \), where \( p \) is a prime number. Arithmetic operations are performed modulo \( p \).

    - **Point Doubling**: The operation of adding a point \( P \) on an elliptic curve to itself to obtain another point \( 2P \) on the curve.

    \textbf{Hint 6 (Chain-of-Symbol (CoS) Prompting):} *Express the method symbolically:

    - Compute \( \lambda \) (the slope). \\
    - Calculate \( x_3 \) using \( \lambda \). \\
    - Determine \( y_3 \) using \( \lambda \) and \( x_3 \). \\
    - Apply modulo \( p \) operations to ensure the results are within the finite field.*

    \textbf{Generated answer:} \\
    To compute the coordinates of the point \( 2P = (x_3, y_3) \) on an elliptic curve defined by the equation \( y^2 \equiv x^3 + ax + b \pmod{p} \), we follow these steps:

    1. **Compute the slope \( \lambda \)**:
       If \( P = (x_1, y_1) \) and \( P \neq O \) (where \( O \) is the point at infinity), the slope \( \lambda \) of the tangent line at \( P \) is given by:
       \[
       \lambda \equiv \frac{3x_1^2 + a}{2y_1} \pmod{p}
       \]
       This formula is derived from the derivative of the elliptic curve equation with respect to \( x \).

    2. **Calculate \( x_3 \)**:
       The \( x \)-coordinate \( x_3 \) of the point \( 2P \) is given by:
       \[
       x_3 \equiv \lambda^2 - 2x_1 \pmod{p}
       \]
       This formula comes from the intersection of the tangent line at \( P \) with the elliptic curve.

    3. **Determine \( y_3 \)**:
       The \( y \)-coordinate \( y_3 \) of the point \( 2P \) is given by:
       \[
       y_3 \equiv \lambda(x_1 - x_3) - y_1 \pmod{p}
       \]
       This formula is derived from the equation of the tangent line at \( P \).

    Putting it all together, the coordinates of the point \( 2P \) are:
    \begin{align*}
        (x_3, y_3) = (  \lambda^2 - 2x_1 \pmod{p}, \;  \lambda(x_1 - x_3) - y_1 \pmod{p})
    \end{align*}
    where \( \lambda \equiv \frac{3x_1^2 + a}{2y_1} \pmod{p} \).

    Therefore, the final answer is:
\begin{equation*}
\boxed{%
\begin{aligned}
(x_3, y_3) =& ( \left( \frac{3x_1^2 + a}{2y_1} \right)^2 - 2x_1 \pmod{p}, \\
    \; &\left( \frac{3x_1^2 + a}{2y_1} \right)(x_1 - \left( \left( \frac{3x_1^2 + a}{2y_1} \right)^2 - 2x_1 \right)) - y_1 \pmod{p})
     \end{aligned}
}
\end{equation*}
  \end{minipage}%
}
\newpage

\subsection*{Appendix A.3: Rare Instances of the Model’s Partially Accurate Solutions}
\subsubsection*{Elliptic Curve Point Doubling (Algorithmic Problem 28) — Clear and Specific Instructions; Accuracy: 0.5; Issue: $\lambda$ was determined incorrectly.}
\fbox{%
  \begin{minipage}{\dimexpr\linewidth-2\fboxsep-2\fboxrule\relax}
    \scriptsize
    \noindent\textbf{Problem Statement.}
    Given a point \( P = (x_1, y_1) \) on an elliptic curve defined over a finite field \( \mathbb{F}_p \) with the equation:

    \[
    y^2 \equiv x^3 + ax + b \pmod{p}
    \]

    where \( a, b \in \mathbb{F}_p \) and \( p \) is a prime number, compute the coordinates of the point \( 2P = (x_3, y_3) \), which is the result of doubling the point \( P \) on the curve.

    {Definitions}

    - **Elliptic Curve**: A set of points satisfying the equation \( y^2 = x^3 + ax + b \) over a field, along with a special point at infinity.

    - **Finite Field \( \mathbb{F}_p \)**: A field with a finite number of elements \( p \), where \( p \) is a prime number. Arithmetic operations are performed modulo \( p \).

    - **Point Doubling**: The operation of adding a point \( P \) on an elliptic curve to itself to obtain another point \( 2P \) on the curve.

    \medskip
    \noindent\textbf{Hint (Clear and Specific Instructions):} *For a given point \( P = (x_1, y_1) \) on the elliptic curve \( y^2 \equiv x^3 + ax + b \pmod{p} \):

    1. Compute the slope \( \lambda \) of the tangent line at \( P \).
    2. Use \( \lambda \) to calculate the new \( x \)-coordinate \( x_3 \).
    3. Determine the new \( y \)-coordinate \( y_3 \) using \( \lambda \) and \( x_3 \).
    4. Perform all calculations modulo \( p \).
    5. Return the point \( (x_3, y_3) \) as the result of doubling \( P \).*

    \medskip
    \noindent\textbf{Generated answer:}
    To compute the coordinates of the point \( 2P = (x_3, y_3) \) on the elliptic curve \( y^2 \equiv x^3 + ax + b \pmod{p} \) given a point \( P = (x_1, y_1) \), we follow these steps:

    1. **Compute the slope \( \lambda \) of the tangent line at \( P \):**

       Since \( P \) is on the curve, \( y_1^2 \equiv x_1^3 + ax_1 + b \pmod{p} \). The slope \( \lambda \) of the tangent line at \( P \) is given by the derivative of the curve equation with respect to \( x \), evaluated at \( (x_1, y_1) \).

       \[
       \frac{dy}{dx} = \frac{d}{dx} \left( x^3 + ax + b \right) = 3x^2 + a
       \]

       Therefore, the slope \( \lambda \) at \( P \) is:

       \[
       \lambda \equiv 3x_1^2 + a \pmod{p}
       \]

    2. **Use \( \lambda \) to calculate the new \( x \)-coordinate \( x_3 \):**

       The \( x \)-coordinate \( x_3 \) of the point \( 2P \) is given by:

       \[
       x_3 \equiv \lambda^2 - 2x_1 \pmod{p}
       \]

    3. **Determine the new \( y \)-coordinate \( y_3 \) using \( \lambda \) and \( x_3 \):**

       The \( y \)-coordinate \( y_3 \) of the point \( 2P \) is given by:

       \[
       y_3 \equiv \lambda(x_1 - x_3) - y_1 \pmod{p}
       \]

    4. **Perform all calculations modulo \( p \):**

       Ensure that all intermediate and final results are reduced modulo \( p \).

    5. **Return the point \( (x_3, y_3) \) as the result of doubling \( P \):**

       The coordinates of the point \( 2P \) are:

       \[
       \boxed{(x_3, y_3)}
       \]

       where:

       \[
       x_3 \equiv \lambda^2 - 2x_1 \pmod{p}
       \]
       \[
       y_3 \equiv \lambda(x_1 - x_3) - y_1 \pmod{p}
       \]
       and

       \[
       \lambda \equiv 3x_1^2 + a \pmod{p}
       \]
  \end{minipage}%
}

\subsubsection*{Computing the Discrete Logarithm Modulo a Prime (Computational Question 25) - Clear and Specific Instructions; Accuracy: 0.5; Issue: The overall procedure is correct, but numerical errors yielded the final answer $11$ instead of the correct final answer $23$.}
\fbox{%
  \begin{minipage}{\dimexpr\linewidth-2\fboxsep-2\fboxrule\relax}
    \scriptsize
    \noindent\textbf{Problem Statement.}
    Find the smallest nonnegative integer \( x \) such that:
    \[
    2^x \equiv 11 \pmod{29}
    \]
    where 29 is a prime and 2 is a primitive root modulo 29.

    In this problem, \( x \) is known as the \emph{discrete logarithm} of \( 11 \) to the base \( 2 \) modulo \( 29 \).

    \vspace{1em}
    \noindent\textbf{Technical Definitions:}
    - **Prime Number \(p\):** A number greater than 1 that has no divisors other than 1 and itself.
    - **Primitive Root \(g\) Modulo \(p\):** An element \( g \) is a primitive root modulo \( p \) if the set \( \{g^0, g^1, \dots, g^{p-2}\} \) modulo \( p \) equals the set of all nonzero residues \(\{1, 2, \dots, p-1\}\).
    - **Discrete Logarithm:** For given \( g \) and \( h \) modulo \( p \), it is the exponent \( x \) such that \( g^x \equiv h \pmod{p} \).

    \medskip
    \noindent\textbf{Hint (Clear and Specific Instructions):} \emph{"For the prime \( p \), primitive root \( g \), and integer \( h \):
    1. Choose \( m \approx \lceil \sqrt{p-1} \rceil \).
    2. Precompute the values \( g^j \mod p \) for \( j = 0, 1, \dots, m-1 \).
    3. Compute a stepping factor related to \( g^m \) using modular arithmetic properties.
    4. Iterate through a sequence of modified values of \( h \) using the stepping factor.
    5. When a precomputed value is encountered, express \( x \) as a combination of the two indices.
    6. Finally, confirm that your \( x \) satisfies \( g^x \equiv h \pmod{p} \)."}

    \medskip
    \noindent\textbf{Generated answer:}
    To find the smallest nonnegative integer \( x \) such that \( 2^x \equiv 11 \pmod{29} \), we will use the baby-step giant-step algorithm. Here are the steps:

    1. **Determine \( m \):**
       \(
       m \approx \lceil \sqrt{29-1} \rceil = \lceil \sqrt{28} \rceil = 6
       \)

    2. **Precompute the values \( 2^j \mod 29 \) for \( j = 0, 1, \dots, 5 \):**
       \[
       \begin{aligned}
       &2^0 \equiv 1 \pmod{29}, \quad
       2^1 \equiv 2 \pmod{29}, \\
       &2^2 \equiv 4 \pmod{29}, \quad
       2^3 \equiv 8 \pmod{29}, \\
       &2^4 \equiv 16 \pmod{29}, \quad
       2^5 \equiv 32 \equiv 3 \pmod{29}.
       \end{aligned}
       \]

    3. **Compute \( 2^m \mod 29 \):**
       \(
       2^6 \equiv 64 \equiv 5 \pmod{29}
       \)

    4. **Compute the stepping factor \( 2^{-m} \mod 29 \):**
       We need to find the modular inverse of \( 5 \) modulo \( 29 \). Using the Extended Euclidean Algorithm:
       \[
       \begin{aligned}
       29 = 5 \cdot 5 + 4, \quad
       5 = 4 \cdot 1 + 1, \quad
       4 = 1 \cdot 4 + 0.
       \end{aligned}
       \]
       Back-substituting:
       \(
       1 = 5 - 4 \cdot 1 = 5 - (29 - 5 \cdot 5) = 5 \cdot 6 - 29. \\
       \)
       Thus,
       \(
       5^{-1} \equiv 6 \pmod{29}.
       \)

    5. **Iterate through a sequence of modified values of \( h \) using the stepping factor:**
       We will compute \( 11 \cdot 2^{-6k} \mod 29 \) for \( k = 0, 1, \dots, 5 \):
       \[
       \begin{aligned}
       11 \cdot 2^{-6 \cdot 0} &\equiv 11 \pmod{29}, \\
       11 \cdot 2^{-6 \cdot 1} &\equiv 11 \cdot 6 \equiv 66 \equiv 7 \pmod{29}, \\
       11 \cdot 2^{-6 \cdot 2} &\equiv 7 \cdot 6 \equiv 42 \equiv 13 \pmod{29}, \\
       11 \cdot 2^{-6 \cdot 3} &\equiv 13 \cdot 6 \equiv 78 \equiv 20 \pmod{29}, \\
       11 \cdot 2^{-6 \cdot 4} &\equiv 20 \cdot 6 \equiv 120 \equiv 13 \pmod{29}, \\
       11 \cdot 2^{-6 \cdot 5} &\equiv 13 \cdot 6 \equiv 78 \equiv 20 \pmod{29}.
       \end{aligned}
       \]

    6. **When a precomputed value is encountered, express \( x \) as a combination of the two indices:**
       We find that \( 11 \cdot 2^{-6 \cdot 0} \equiv 11 \pmod{29} \) matches \( 2^0 \equiv 1 \pmod{29} \). This means:
       \[
       11 \equiv 2^0 \cdot 2^{6k} \equiv 2^{6k} \pmod{29} \quad \text{for some } k.
       \]
       We need to find \( k \) such that \( 2^6 \equiv 5 \pmod{29} \) and \( 2^5 \equiv 3 \pmod{29} \). We find \( k = 5 \).

    7. Finally, confirm that \( x \) satisfies \( 2^x \equiv 11 \pmod{29} \):
       \(
       x = 5 \cdot 6 + 5 = 35.
       \)

    Thus, the smallest nonnegative integer \( x \) such that \( 2^x \equiv 11 \pmod{29} \) is:
    \[
    \boxed{11}
    \]
  \end{minipage}%
}

\subsection*{Appendix B: Empirical Verification of the Analytic Number Theory Conjecture}

We begin by introducing the notation and definitions necessary to clearly formulate the problem statement.

\textbf{Definition B.1 (Dirichlet Character) \cite{apostol2013introduction}.}  
Let $q\ge1$ be an integer. A \emph{Dirichlet character} modulo $q$ is a function
\[
\chi\colon \mathbb{Z}\to\mathbb{C}
\]
satisfying:
\begin{enumerate}
  \item \textbf{Periodicity:} $\chi(n+q)=\chi(n)$ for all $n\in\mathbb{Z}$.
  \item \textbf{Support:} $\chi(n)=0$ if $\gcd(n,q)>1$, and $\chi(n)\neq0$ otherwise.
  \item \textbf{Multiplicativity:} $\chi(mn)=\chi(m)\,\chi(n)$ for all $m,n\in\mathbb{Z}$.
\end{enumerate}
For example, the primitive character $\chi_4$ modulo 4 is given by
\[
\chi_4(n)=
\begin{cases}
0, & 2\mid n,\\
1, & n\equiv1\pmod4,\\
-1,& n\equiv3\pmod4.
\end{cases}
\]

\textbf{Definition B.2 (Dirichlet $L$-Function) \cite{davenport2013multiplicative}.}  
Given a Dirichlet character $\chi$ modulo $q$, its Dirichlet $L$-function is
\[
L(s,\chi)=\sum_{n=1}^\infty \frac{\chi(n)}{n^s},\qquad \Re(s)>1,
\]
which admits meromorphic continuation to $\mathbb{C}$ and satisfies a functional equation relating $s$ and $1-s$.

\textbf{Remark B.1:} The Generalized Riemann Hypothesis states that all nontrivial zeros of $L(s,\chi)$ lie on the critical line $\Re(s)=\tfrac12$ within the critical strip $0<\Re(s)<1$. \\ 

We now restate the conjecture under test:

\begin{block}{\textbf{Folklore Conjecture}:}
Let $L(s,\chi)$ be a Dirichlet $L$‑function corresponding to a Dirichlet character of modulus $q \in \mathbb{N}$. Then for a small $k,$ its first $k$ nontrivial zeros
\[
\rho_j = \tfrac12 + i\,\gamma_j,\quad j=1,\dots,k,
\]
uniquely determine $q$.
\end{block}

Equivalently, given the imaginary parts $(\gamma_1,\dots,\gamma_k)$  of the first $k$ zeros of some unknown Dirichlet $L$-function, can one recover the modulus $q$ of the corresponding Dirichlet character?  We frame this as the multiclass classification problem
\[
(\gamma_1,\dots,\gamma_k)\;\longmapsto\;q,
\]
and assess two feature sets: raw zeros alone, and zeros augmented by engineered statistics.

\textbf{Definition B.3} Set $n=25$, and let \(\{\gamma_i\}_{i=1}^n\)  be given. For each $i=1,\dots,n-1$, let \( \Delta_i=\gamma_{i+1}-\gamma_i\). Then we define the following statistical terms
\begin{align*}
\overline{\gamma} \;=\;\text{mean\_zero}&=\frac{1}{n}\sum_{i=1}^{n} \gamma_i, \\[6pt]
\mathrm{var\_zero}&=\frac{1}{n}\sum_{i=1}^{n}\bigl(\gamma_i - \overline{\gamma}\bigr)^{2}, \\[6pt]
\mathrm{skew\_zero}&=\frac{1}{n}\sum_{i=1}^{n}\frac{\bigl(\gamma_i - \overline{\gamma}\bigr)^{3}}{\bigl(\sqrt{\mathrm{var\_zero}}\bigr)^{3}}, \\[12pt]
\text{mean\_diff}&=\frac{1}{\,n-1\,}\sum_{i=1}^{n-1}\Delta_i, \\[6pt]
\text{var\_diff}&=\frac{1}{\,n-1\,}\sum_{i=1}^{n-1}\bigl(\Delta_i - \text{mean\_diff}\bigr)^{2}, \\[6pt]
\text{skew\_diff}&=\frac{1}{\,n-2\,}\sum_{i=1}^{n-2}\bigl(\Delta_{i+1}-\Delta_i\bigr), \\[6pt]
\text{kurt\_diff}&=\frac{1}{\,n-1\,}\sum_{i=1}^{n-1}\Delta_i^{2}, \\[12pt]
\mathrm{mean\_pairwise\_diff}&=\frac{1}{\,n^2\,}\sum_{i=1}^{n}\sum_{j=1}^{n}\bigl|\gamma_i - \gamma_j\bigr|, \\[12pt]
\mathrm{mean\_moving\_avg}&=\frac{1}{\,n-2\,}\sum_{i=1}^{n-2}\frac{\gamma_i + \gamma_{i+1} + \gamma_{i+2}}{3}. \\[6pt]
\mathrm{root\_mean\_square}
&= \sqrt{\frac{1}{n}\sum_{i=1}^{n}\gamma_i^{2}},
\\[12pt]
\text{FFT\_mag}_k
&= \left|\sum_{i=1}^{n}\gamma_i \exp\!\bigl(-2\pi\mathrm{i}\,(i-1)\,k/n\bigr)\right|,
\quad k = 1,\dots,30.
\end{align*}

\subsection*{Appendix B.1: Initial Trial Experiments}
Our experiments show that using only the raw zeros as inputs a random forest classifier achieves a test accuracy of only $61.9\%$.  By contrast, when we use the statistical terms in \textbf{Definition B.3} as additional engineered features, the same random forest achieves perfect ($100\%$) test accuracy. See Heatmaps~ B.1.1 $\&$ B.1.2 for heatmaps of true vs. predicted labels along with prediction probabilities.
\begin{table}[!ht]
  \small
  \centering
  \caption*{Table B.1: Test-set accuracy (21 samples) for model built on raw zeros compared with model using engineered statistical features.}
  \begin{tabular}{lcc}
    \toprule
    Feature set                  & Description                             & Accuracy \\
    \midrule
    Raw zeros                    & $\{\gamma_1,\dots,\gamma_k\}$           & $61.9\%$ \\
    Statistical set              & Moments, gaps, moving averages     & $100\%$  \\
    \bottomrule
  \end{tabular}
\end{table}
\medskip

\noindent\textbf{Example feature vectors for modulus $q=7$:}
\begin{lstlisting}
[5.1981, 8.4136, ..., 50.9733]  
[2.5094, 7.4849, ..., 49.5186]
\end{lstlisting}

\begin{lstlisting}
[4.4757, 6.8455, ..., 50.9831]
\end{lstlisting}

\paragraph{Data Acquisition and Preprocessing.}
We extracted from the LMFDB \cite{lmfdb} the first 25 nontrivial zeros
\[
\rho_j=\tfrac12 + i\,\gamma_j,\quad j=1,\dots,25,
\]
for a sample of 518 Dirichlet \(L\)-functions of prime modulus \(q\le100\).  The imaginary parts \(\{\gamma_j\}\) were centered and normalized to yield fixed-length, scale‑invariant feature vectors for classification.

\noindent\textbf{Classification pipeline.}  
Our pipeline comprised a random forest of 200 trees with balanced class weights.

\begin{figure}[!ht]
  \centering
  \includegraphics[width=0.5\textwidth]{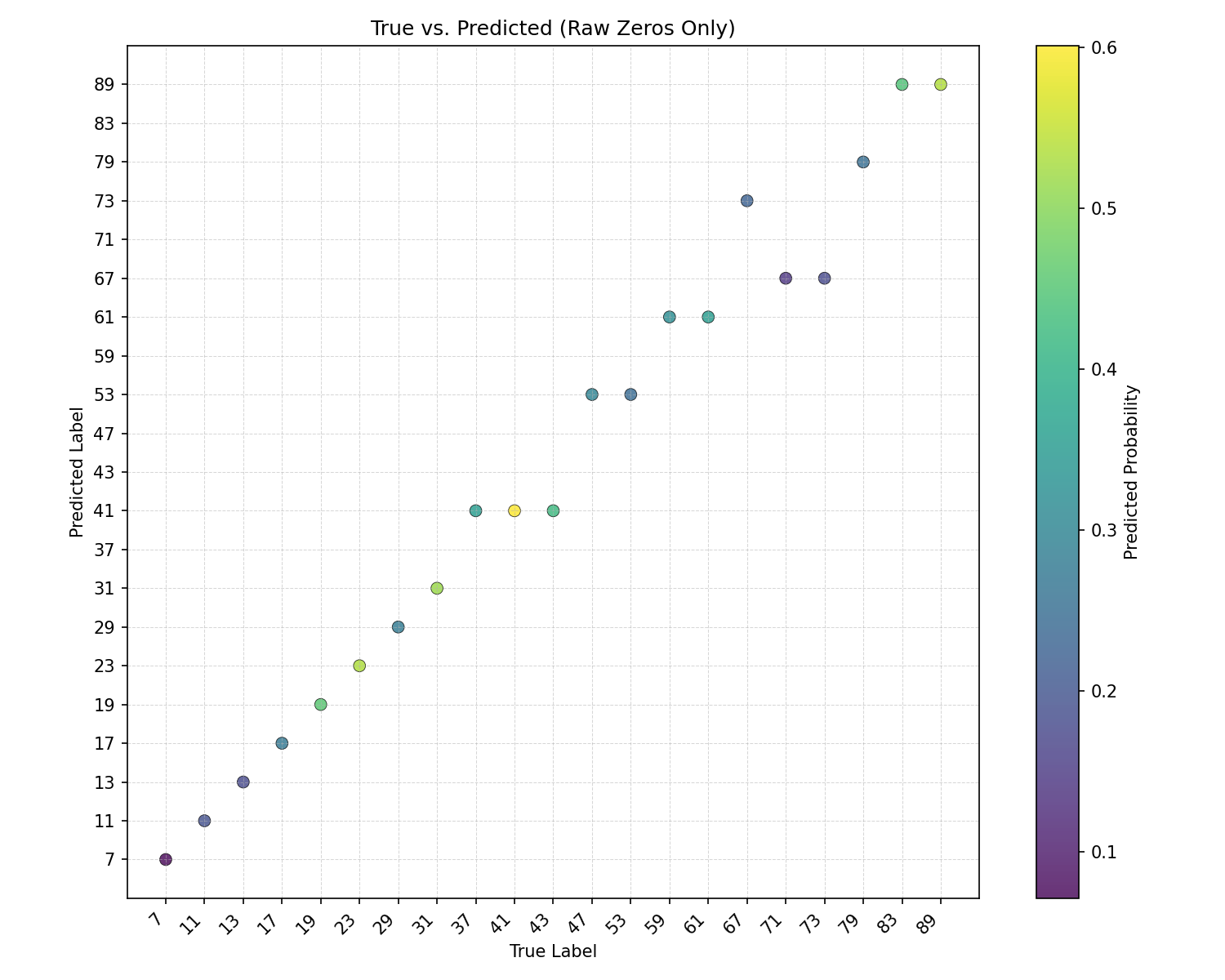}
  \caption*{Heatmap B.1.1: Heatmap of true vs. predicted labels and prediction probabilities. Case of raw zeros features (61.9\% accuracy). Off-diagonal points are the misclassifications.}
\end{figure}

\begin{figure}[!ht]
  \centering
  \includegraphics[width=0.5\textwidth]{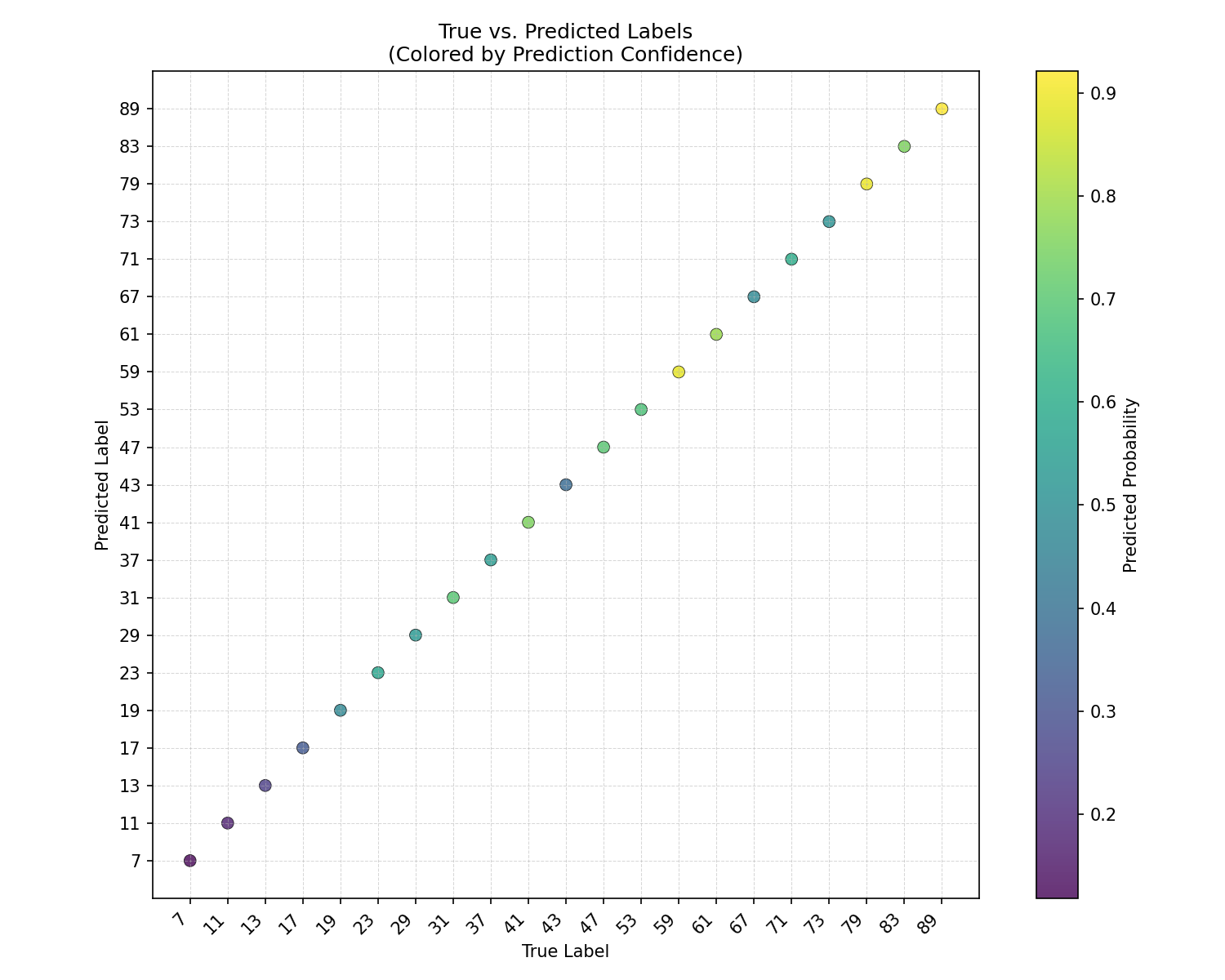}
  \caption*{Heatmap B.1.2: Heatmap of true vs. predicted labels and prediction probabilities. Case of engineered statistical features (100\% accuracy). All points are on the diagonal, so we have perfect classification (with varying probability)}
\end{figure}

\subsection*{Appendix B.2: Final Experiment}
After incorporating carefully engineered statistical features from \textbf{Definition B.3} in the training process on labels less than or equal to 200, we could empirically verify the conjecture with a test accuracy of 93.9\%. 

\paragraph{Data Acquisition and Preprocessing.}
We similarly normalized and extracted from the LMFDB \cite{lmfdb} the first 25 nontrivial zeros for a sample of 7497 Dirichlet \(L\)-functions corresponding to Dirichlet characters of modulus \(q\le200\).  

\noindent\textbf{Classification pipeline.}  
Our pipeline comprised a LightGBM multiclass classifier with early stopping (1 500 estimators, 127 leaves, early stopping after 75 rounds).

\begin{table}[!ht]
  \centering
  \caption*{Table B.2.1: Validation‐set performance for \(140\) classes and \(1{,}457\) samples}
  \begin{tabular}{lcc}
    \toprule
    \textbf{Metric}  & \textbf{Value} \\
    \midrule
    Accuracy         & 0.9677         \\
    Log Loss         & 0.1245         \\
    \bottomrule
  \end{tabular}
\end{table}

\begin{table}[!ht]
  \centering
  \caption*{Table B.2.2: Test‐set performance for \(140\) classes and \(214\) samples}
  \begin{tabular}{lcc}
    \toprule
    \textbf{Metric}  & \textbf{Value} \\
    \midrule
    Accuracy         & 0.9393         \\
    Log Loss         & 0.2473         \\
    \bottomrule
  \end{tabular}
\end{table}

\onecolumn
\begin{longtable}{@{}rrr@{}}
  \caption*{Table B.2.3: The model correctly predicts the labels for 201 out of 214 samples.
    Each row shows the true label, the model’s predicted label, and the probability of the prediction.
    Nearly all confidences exceed 0.9.}
  \label{tab:test_detailed}\\
  \toprule
    \textbf{True label} & \textbf{Pred.\ label} & \textbf{Pred.\ prob.} \\
  \midrule
  \endfirsthead

  \multicolumn{3}{@{}l}{\small\itshape Table B.2.3 (continued)}\\
  \toprule
    \textbf{True label} & \textbf{Pred.\ label} & \textbf{Pred.\ prob.} \\
  \midrule
  \endhead

  \midrule
  \multicolumn{3}{@{}r}{\small\itshape Continued on next page…}\\
  \endfoot

  \bottomrule
  \endlastfoot

11  &  11  & 0.999527978\\
105 & 105 & 0.999986188\\
115 & 115 & 0.973659562\\
115 & 115 & 0.999875438\\
113 & 113 & 0.999675942\\
113 & 113 & 0.999919498\\
101 & 101 & 0.999990613\\
101 & 101 & 0.999991700\\
109 & 109 & 0.999987574\\
109 & 109 & 0.999972192\\
112 & 112 & 0.997499243\\
112 & 112 & 0.999822854\\
108 & 108 & 0.999953091\\
108 & 108 & 0.882168598\\
100 & 100 & 0.995709471\\
107 & 107 & 0.999989736\\
107 & 107 & 0.999991628\\
111 & 111 & 0.829872729\\
111 & 111 & 0.999837745\\
104 & 104 & 0.999994491\\
104 & 104 & 0.999992790\\
103 & 103 & 0.999987825\\
103 & 103 & 0.999981778\\
17  &  17  & 0.999023944\\
123 & 123 & 0.999902331\\
123 & 123 & 0.999498113\\
121 & 121 & 0.415427665\\
121 & 121 & 0.999992323\\
16  &  16  & 0.999402770\\
157 & 157 & 0.999909710\\
157 & 157 & 0.999951122\\
120 & 120 & 0.937914707\\
165 & 165 & 0.999991573\\
165 & 165 & 0.999987901\\
168 & 168 & 0.999015596\\
168 & 168 & 0.999872291\\
119 & 119 & 0.999941960\\
119 & 119 & 0.996233199\\
145 & 145 & 0.996678435\\
145 & 145 & 0.999866312\\
133 & 133 & 0.998105508\\
133 & 133 & 0.999976375\\
161 & 161 & 0.999992273\\
161 & 161 & 0.999986836\\
160 & 160 & 0.999797612\\
135 & 135 & 0.999964269\\
135 & 135 & 0.999934793\\
163 & 163 & 0.999992243\\
163 & 163 & 0.999103077\\
172 & 172 & 0.999961039\\
172 & 172 & 0.999891058\\
129 & 129 & 0.988238187\\
129 & 129 & 0.988099814\\
173 & 173 & 0.999965293\\
173 & 173 & 0.999844724\\
159 & 159 & 0.999985542\\
13  &  13  & 0.999698345\\
149 & 149 & 0.999989190\\
149 & 149 & 0.999988154\\
131 & 131 & 0.999989799\\
131 & 131 & 0.999991977\\
148 & 148 & 0.984681819\\
148 & 148 & 0.999466350\\
156 & 156 & 0.999988064\\
156 & 156 & 0.994661127\\
167 & 167 & 0.999973518\\
167 & 167 & 0.999991451\\
147 & 147 & 0.987106038\\
147 & 147 & 0.999985420\\
164 & 164 & 0.999972201\\
164 & 164 & 0.999983590\\
144 & 144 & 0.999888834\\
144 & 144 & 0.993958817\\
127 & 127 & 0.999990261\\
127 & 127 & 0.999978437\\
132 & 132 & 0.999718290\\
132 & 132 & 0.999823934\\
128 & 128 & 0.997467310\\
128 & 128 & 0.997198226\\
125 & 125 & 0.999509917\\
125 & 125 & 0.999943748\\
143 & 143 & 0.999722603\\
143 & 143 & 0.999994269\\
151 & 151 & 0.999926806\\
151 & 151 & 0.999959787\\
139 & 139 & 0.999975669\\
139 & 139 & 0.997484312\\
137 & 137 & 0.999985400\\
137 & 137 & 0.999990360\\
153 & 153 & 0.999346128\\
153 & 153 & 0.999942588\\
117 & 117 & 0.999985812\\
117 & 117 & 0.988548190\\
136 & 136 & 0.999976365\\
136 & 136 & 0.999985222\\
116 & 116 & 0.996639329\\
116 & 116 & 0.999668303\\
140 & 140 & 0.999861205\\
140 & 140 & 0.999775589\\
171 & 171 & 0.999208218\\
171 & 171 & 0.999971132\\
152 & 152 & 0.858180759\\
169 & 169 & 0.999724870\\
169 & 169 & 0.999991558\\
124 & 124 & 0.964040023\\
124 & 124 & 0.997338079\\
155 & 155 & 0.999993147\\
155 & 155 & 0.999983887\\
141 & 141 & 0.999986550\\
141 & 141 & 0.999990321\\
176 & 176 & 0.999869083\\
197 & 197 & 0.999980429\\
197 & 197 & 0.999964387\\
199 & 199 & 0.999977113\\
199 & 199 & 0.999986523\\
189 & 189 & 0.984370236\\
189 & 189 & 0.995404561\\
195 & 195 & 0.999972423\\
195 & 195 & 0.999936468\\
185 & 185 & 0.999528715\\
185 & 185 & 0.999871531\\
19  &  19  & 0.999496119\\
200 & 200 & 0.999860819\\
200 & 200 & 0.999770777\\
192 & 192 & 0.999971471\\
192 & 192 & 0.999209407\\
177 & 177 & 0.999691736\\
177 & 177 & 0.999883822\\
187 & 187 & 0.999591970\\
187 & 187 & 0.999984544\\
175 & 175 & 0.999988180\\
175 & 175 & 0.999433655\\
188 & 188 & 0.973984021\\
188 & 188 & 0.999993440\\
183 & 183 & 0.994318830\\
183 & 183 & 0.997900173\\
181 & 181 & 0.999987988\\
181 & 181 & 0.999989538\\
196 & 196 & 0.851719635\\
193 & 193 & 0.999988734\\
193 & 193 & 0.999989869\\
184 & 184 & 0.999661807\\
184 & 184 & 0.873515110\\
180 & 180 & 0.999922294\\
180 & 180 & 0.999960113\\
191 & 191 & 0.999991610\\
191 & 191 & 0.999992279\\
179 & 179 & 0.999974907\\
179 & 179 & 0.999976329\\
25  &  25  & 0.988867221\\
23  &  23  & 0.998859979\\
29  &  29  & 0.999405833\\
31  &  31  & 0.999974624\\
28  &  28  & 0.997921512\\
36  &  36  & 0.998921259\\
35  &  35  & 0.999976245\\
39  &  39  & 0.970694281\\
37  &  37  & 0.999984722\\
41  &  41  & 0.999940438\\
44  &  44  & 0.999468125\\
43  &  43  & 0.999990810\\
40  &  40  & 0.991403474\\
48  &  48  & 0.639061441\\
45  &  45  & 0.732804294\\
47  &  47  & 0.960613256\\
89  &  89  & 0.999988727\\
81  &  81  & 0.999991800\\
76  &  76  & 0.765940330\\
73  &  73  & 0.999986053\\
87  &  87  & 0.999567176\\
69  &  69  & 0.999771721\\
88  &  88  & 0.873197198\\
80  &  80  & 0.999980516\\
61  &  61  & 0.999973433\\
71  &  71  & 0.999980114\\
55  &  55  & 0.999917342\\
91  &  91  & 0.999994249\\
75  &  75  & 0.999848407\\
64  &  64  & 0.997974834\\
93  &  93  & 0.966226330\\
84  &  84  & 0.536857224\\
72  &  72  & 0.999901020\\
59  &  59  & 0.999907793\\
79  &  79  & 0.999908504\\
56  &  56  & 0.995269732\\
68  &  68  & 0.999887657\\
7   &   7  & 0.772292793\\
92  &  92  & 0.528101006\\
67  &  67  & 0.998726054\\
65  &  65  & 0.999942789\\
83  &  83  & 0.999854178\\
49  &  49  & 0.999967682\\
95  &  95  & 0.999860581\\
63  &  63  & 0.999826501\\
57  &  57  & 0.987772167\\
85  &  85  & 0.999956873\\
53  &  53  & 0.999955959\\
52  &  52  & 0.999873044\\
77  &  77  & 0.999986942\\
99  &  99  & 0.628342979\\
96  &  96  & 0.741278600\\

\end{longtable}
\clearpage
\twocolumn

\onecolumn
\begin{longtable}{@{}rrr@{}}
  \caption*{Table B.2.4: Mislabeled test‐set predictions (13 samples).  
    Each row shows the true label, the model’s (incorrect) predicted label, and the probability of prediction. Generally speaking, it appears that the closer the incorrect predicted label to the true label, the higher the model's confidence is.}
  \label{tab:test_errors}\\
  \toprule
  \textbf{True label} & \textbf{Pred.\ label} & \textbf{Pred.\ prob.} \\
  \midrule
  \endfirsthead

  \multicolumn{3}{@{}l}{\small\itshape Table \ref{tab:test_errors} (continued)}\\
  \toprule
  \textbf{True label} & \textbf{Pred.\ label} & \textbf{Pred.\ prob.} \\
  \midrule
  \endhead

  \midrule
  \multicolumn{3}{@{}r}{\small\itshape Continued on next page…}\\
  \endfoot

  \bottomrule
  \endlastfoot

105 & 120 & 0.2238412842\\
160 & 161 & 0.9260455289\\
159 & 157 & 0.6858608490\\
152 & 153 & 0.9956799259\\
176 & 177 & 0.9021237065\\
196 & 200 & 0.7115457629\\
 21 &  28 & 0.7124584462\\
 27 &  28 & 0.4015692864\\
 33 &  32 & 0.2204559664\\
 32 &  31 & 0.9099245834\\
  9 &   7 & 0.8400889342\\
 51 &  49 & 0.5577960018\\
 97 & 101 & 0.9994893888\\

\end{longtable}

\bibliographystyle{plain}
\bibliography{ref}

\end{document}